\documentclass[amsthm]{elsart}
\usepackage{amsfonts,amssymb,amsmath}%,showlabels}

\def\diam{\mbox{\rm diam\,}}
\newtheorem{theorem}[thm]{Theorem}
\newtheorem{corollary}[thm]{Corollary}
\newtheorem{example}[thm]{Example}
\newtheorem{question}[thm]{Question}
\newtheorem{lemma}[thm]{Lemma}

\newtheorem{definition}[thm]{Definition}
\newtheorem{remark}[thm]{Remark}
\newtheorem{remarks}[thm]{Remarks}
\date{}
\begin{document}
\begin{frontmatter}
\title{The convergence space of minimal usco mappings}

\author{R. Anguelov\thanksref{NRF}}
\address{University of Pretoria, Department of Mathematics and Applied Mathematics\\ Pretoria 0002, South Africa}
\ead{roumen.anguelov@up.ac.za}
 and
\author{O.F.K. Kalenda\thanksref{OK}}
\address{Charles University, Faculty of Mathematics and Physics,
 Department of Mathematical Analysis, Sokolovsk\'a~83,
186~75, Praha~8, Czech Republic} \ead{kalenda@karlin.mff.cuni.cz}

\thanks[NRF]{Supported in part by Research grant NRF4886.}

\thanks[OK]{Supported in part by Research project MSM 0021620839 financed by
MSMT, by Research grant GA\v{C}R 201/06/0018 and by Universidad
Polit\'ecnica de Valencia.}
\begin{abstract}
A convergence structure generalizing the order convergence structure
on the set of Hausdorff continuous interval functions is defined on
the set of minimal usco maps. The properties of the obtained
convergence space are investigated and essential links with the
pointwise convergence and the order convergence are revealed. The
convergence structure can be extended to a uniform convergence
structure so that the convergence space is complete. The important
issue of the denseness of the subset of all continuous functions is
also addressed.
\end{abstract}

\end{frontmatter}

\section{Introduction}

The aim of this paper is to define a convergence structure on the
set $\mathcal{M}(X,Y)$ of minimal usco maps from from $X$ to $Y$
where $X$ and $Y$ are given topological spaces and to establish
basic properties of the obtained convergence space. For
$Y=\mathbb{R}$, the minimal usco maps can be identified with the
Hausdorff continuous interval valued functions on $X$. Therefore the
set $\mathcal{M}(X,Y)$ can be considered as a generalization of the
set Hausdorff continuous functions for arbitrary topological spaces
$Y$. Hence we consider on $\mathcal{M}(X,Y)$ a convergence structure
which generalizes the order convergence structure on the set of
Hausdorff continuous functions, see \cite{AnguelovWalt}. The origin
of the Hausdorff continuous functions is in the theory of Hausdorff
approximation. However, they have a particular significance in the
general theory of PDEs. It was shown recently, see \cite{AngRos1},
\cite{AngRos2}, that the solutions of large classes of nonlinear
PDEs can be assimilated with Hausdorff continuous functions. It is
expected that the convergence space $\mathcal{M}(X,Y)$ derived and
analyzed in this paper will further facilitate this development
particularly for multidimensional problems.

The paper is organized as follows.

For completeness of the exposition we recall in the next section the
definitions of usco and minimal usco maps and give some of their
basic properties. We also define the notion of a quasiminimal usco
map which is  quite important for the topic.

In Section \ref{sec-def} we define convergence of filters on the
space $\mathcal{M}(X,Y)$ and we prove it satisfies the axioms of a
convergence structure. Some basic properties together with
characterizations of the convergent sequences and nets are also
presented. It is shown that the convergence is not topological.
Hence $\mathcal{M}(X,Y)$ is a convergence space but not a
topological space.

The relationship of the convergence in $\mathcal{M}(X,Y)$ and the
pointwise convergence is studied in Section \ref{sec-pointwise}. It
is shown through examples that in general neither convergence
implies the other. Nevertheless a strong connection exists. In
particular, for $X$ a Baire space and $Y$ a metric space any filter
convergent in $\mathcal{M}(X,Y)$ converges pointwise on a residual
subset of $X$.

In Section \ref{sec-order} we consider the special case when the
target space $Y$ is the real line. Then $\mathcal{M}(X,\mathbb{R})$
can be ordered similarly to the way interval functions are ordered
and we show that the convergence in $\mathcal{M}(X,\mathbb{R})$ is
equivalent to the order convergence. Hence
$\mathcal{M}(X,\mathbb{R})$ is isomorphic to the convergence space
of Hausdorff continuous functions on $X$ equipped with the order
convergence structure.

Section \ref{sec-unif} contains the definition of a uniform
convergence structure on $\mathcal{M}(X,Y)$ for the case when $X$ is
a Baire space and $Y$ a metric space. We show that this uniform
convergence structure induces our convergence structure and that
$\mathcal{M}(X,Y)$ is complete.

In Section \ref{sec-dens} we consider the set $\mathcal{C}(X,Y)$ of
all continuous functions. The concept of minimal usco generalizes
the concept of continuity while retaining some of its essential
properties. It is interesting from both theoretical and practical
points of view when the set $\mathcal{C}(X,Y)$ is dense in
$\mathcal{M}(X,Y)$. We give a partial answer formulating some open
questions as well.

\section{Usco and minimal usco maps}\label{sec-usco}

Let $X$ and $Y$ be topological spaces. A set-valued map $g:X\to Y$ is called
{\it upper semicontinuous compact valued\/} (shortly \textit{usco}) if
\begin{itemize}
    \item $g(x)$ is a nonempty compact subset of $Y$ for each $x\in X$;
    \item $\{ x\in X : g(x)\subset U\}$ is open in $X$ for each open subset $U$ of $Y$.
\end{itemize}

We will always assume that the range space $Y$ is Hausdorff. For the domain space $X$ we require no separation axioms.

A set-valued map $g:X\to Y$ is canonically identified with its {\it graph}, i.e. with the set
$$\{(x,y)\in X\times Y : y\in g(x) \}.$$
Using this identification we will consider unions, intersections and inclusions of set-valued mappings. Hence, for example, if $g:X\to Y$ and $h:X\to Y$ are two set-valued mappings, then $g\subset h$ means that the graph of $g$ is a subset of the graph of $h$, i.e., $g(x)\subset h(x)$ for each $x\in X$.

If $g:X\to Y$ is a set-valued mapping and $A\subset X$ we use,
following the standard convention, the symbol $g(A)$ to denote
$\bigcup\{g(x):x\in A\}$.

We will need the following basic stability properties of usco maps.

\begin{lemma}\label{usco-stability}\quad
\begin{itemize}
    \item[(i)] Let $g_j:X \to Y$, $j=1,\dots,n$, be usco maps. Then $g_1\cup\dots\cup g_n$ is usco as well.
    \item[(ii)] Let $\mathcal{G}$ be a family of usco maps from $X$ to $Y$ such that for each finite subfamily $\mathcal{K}\subset\mathcal{G}$ the intersection $\bigcap\mathcal{K}$ is a nonempty-valued mapping. Then $\bigcap\mathcal{G}$ is usco.
\end{itemize}
\end{lemma}

\begin{proof} (i) Let $g=g_1\cup\dots\cup g_n$. Then $g(x)=g_1(x)\cup\dots\cup g_n(x)$ for each $x\in X$ and hence it is a nonempty compact set. Further, if $x\in X$ and $U\subset Y$ is open such that $g(x)\subset U$, then $g_j(x)\subset U$ for each $j\in\{1,\dots,n\}$. Hence there are $V_j$, $j\in\{1,\dots,n\}$ neighborhoods of $x$ in $X$ such that $g_j(V_j)\subset U$ for each $j$. Then $V=V_1\cap\dots\cap V_n$ is a neighborhood of $x$ in $X$ satisfying $g(V)\subset U$.

(ii) First suppose that $g_1$ and $g_2$ are usco maps and $g=g_1\cap g_2$ is  nonempty-valued.
Then $g(x)=g_1(x)\cap g_2(x)$ is a nonempty compact set for each $x\in X$. Let $x\in X$ and $U$ be an open subset of $Y$ such that $g(x)\subset U$. Then $g_1(x)\setminus U$ and $g_2(x)\setminus U$ are disjoint compact subsets of $Y$. As $Y$ is Hausdorff, there are disjoint open sets $U_1$ and $U_2$ such that
$$g_1(x)\setminus U\subset U_1 \qquad\mbox{and}\qquad g_2(x)\setminus U\subset U_2.$$
Then for each $j\in\{1,2\}$ we have $g_j(x)\subset U\cup U_j$ and hence there is a neighborhood $V_j$ of $x$ with $g_j(V_j)\subset U\cup U_j$. Set $V=V_1\cap V_2$. Then $V$ is a neighborhood of $x$ and $g(V)\subset (U\cup U_1)\cap (U\cup U_2)=U$.
Thus $g$ is usco.

We have just proved the assertion (ii) for two-element families $\mathcal{G}$. It is easy to show by induction that it is true for any finite $\mathcal G$.

Suppose now that $\mathcal G$ is an arbitrary family of usco mappings such that each finite subfamily has nonempty-valued intersection. Set $h=\bigcap\mathcal G$. As for each $x\in X$ we have $h(x)=\bigcap_{g\in \mathcal{G}} g(x)$ and the family $\{g(x):g\in \mathcal{G}\}$ is by our assumptions centered, $h(x)$ is a nonempty compact subset of $Y$.

Further, let $x\in X$ and $U$ be an open subset of $Y$ such that $h(x)\subset U$.
Then there is a finite subfamily $\mathcal{K}\subset\mathcal{G}$ such that
$$\bigcap_{g\in\mathcal{K}} g(x)\subset U.$$
The map $\tilde g=\bigcap\mathcal{K}$ is usco by the above and satisfies $\tilde g(x)\subset U$. Hence there is a neighborhood $V$ of $x$ such that $\tilde g(V)\subset U$. Then clearly $h(V)\subset U$ and the proof is finished.
\end{proof}

We will also need the following lemma on modifying usco maps. The proof is trivial and hence we omit it.

\begin{lemma}\label{usco-modif}
Let $g:X\to Y$ be a usco map.
\begin{itemize}
    \item[(i)] Let $h\subset g$ be a set-valued mapping. Suppose there is an open set $U\subset X$ such that
      \begin{itemize}
                \item[$\bullet$] $h(x)=g(x)$ for $x\in X\setminus U$;
                \item[$\bullet$] $h|_U:U\to Y$ is usco.
            \end{itemize}
        Then $h$ is usco.
   \item[(ii)] Let $U\subset X$ be open and $F\subset Y$ be closed. Then the mapping
          $h:X\to Y$ defined by
          $$ h(x)=\begin{cases} g(x)\cap F, & x\in U,\\ g(x), & x\in X\setminus U, \end{cases}$$
          is usco provided it is nonempty-valued.
\end{itemize}
\end{lemma}

A usco map $g:X\to Y$ is called {\it minimal} if it is minimal with
respect to inclusion, i.e., if $g=h$ whenever $h:X\to Y$ is usco
satisfying $h\subset g$. It is a well-know consequence of Zorn's
lemma (and of Lemma~\ref{usco-stability}(ii)) that for each usco map
$g:X\to Y$ there is a minimal usco $h\subset g$,
\cite{BorweinKortezov}.

The following characterization of minimal usco maps is proved in \cite[Lemma 3.1.2]{fabian}

\begin{lemma}\label{minimal-prop}
Let $g:X\to Y$ be a usco map. The following assertions are equivalent.
\begin{itemize}
    \item[(i)] $g$ is minimal.
    \item[(ii)] Whenever $V\subset X$ and $U\subset Y$ are open sets such that $g(V)\cap U\ne\emptyset$, there is a nonempty open set $W\subset V$ with $g(W)\subset U$.
\end{itemize}
\end{lemma}
We will denote by $\mathcal{M}(X,Y)$ the set of all minimal usco
maps from $X$ to $Y$. The minimal usco maps generalize the concept
of continuous function and retain some of its properties. For
example, a minimal usco map is completely determined by its values
on a dense subset of the domain as stated in the next lemma.

\begin{lemma}\label{tfgequal}
Let $f$ and $g$ be usco mappings from $X$ to $Y$ such that $f$ is
minimal. If $f\not\subset g$, then there is a nonempty open set
$U\subset X$ such that $f(U)\cap g(U)=\emptyset$.

In particular, if  $f,g\in\mathcal{M}(X,Y)$ are such that there
exists a dense subset $D$ of $X$ such that $f(x)\cap
g(x)\ne\emptyset$ for each $x\in D$, then $f=g$.
\end{lemma}

\begin{proof}
Suppose that $f$ is minimal and $f\not\subset g$. If $f(x)\cap
g(x)\neq\emptyset$ for all $x\in X$ then $f\cap g$ is (by
Lemma~\ref{usco-stability}) an usco map contained in both $f$ and
$g$. Hence $f=f\cap g\subset g$. Therefore, there exists $x_0\in X$
such that $f(x_0)\cap g(x_0)=\emptyset$. Since $Y$ is Hausdorff
there exist disjoint open subsets $V_1$ and $V_2$ of $Y$ such that
$f(x_0)\subset V_1$ and $g(x_0)\subset V_2$. Using that $f$ and $g$
are usco maps there exists an open neighborhood $U$ of $x_0$ such
that $f(U)\subseteq V_1$ and $g(U)\subseteq V_2$. Hence $f(U)\cap
g(U)=\emptyset$.

Now suppose that $f$ and $g$ are minimal and that $f(x)\cap g(x)$ is
nonempty for all $x$ from a dense subset of $X$. By the first part
of the lemma we get $f\subset g$ and $g\subset f$, hence $f=g$.
\end{proof}

Let $g$ be an usco map from $X$ to $Y$. We associate with $g$ the
following subset of $\mathcal{M}(X,Y)$:
\begin{equation}\label{[g]}
[g]=\{f\in \mathcal{M}(X,Y):f\subset g\}
\end{equation}
By the above the set $[g]$ is not empty.
If $g$ is a minimal usco map then we have $[g]=\{g\}$.
On the other hand, if $[g]$ contains only one element, it need not be minimal.
Since such usco maps will be important for us, we call them {\it quasiminimal}.

The following lemmata present some properties of quasiminimal usco
maps which we will need in the sequel.

\begin{lemma}\label{tQminDense}
Let $g$ be an usco map from $X$ to $Y$. If there exists a dense
subset $D$ of $X$ such that $g(x)$ is a singleton for all $x\in D$
then $g$ is quasiminimal.
\end{lemma}

\begin{proof}Let $f_1,f_2\in[g]$. Since $g$ is singlevalued on $D$
we have $f_1(x)=f_2(x)=g(x)$, $x\in D$, and by Lemma
\ref{tfgequal} we obtain $f_1=f_2$.
\end{proof}

\begin{lemma}\label{union}
Let $g_1$ and $g_2$ be quasiminimal usco maps with $[g_1]=[g_2]$. Then $g_1\cup g_2$ is quasiminimal, too. (And $[g_1\cup g_2]=[g_1]$.)
\end{lemma}

\begin{proof}
The map $g_1\cup g_2$ is usco by Lemma~\ref{usco-stability}.
Let $f$ be the unique element of $[g_1]$. Then clearly $f\in [g_1\cup g_2]$.
We will show it is the unique element.

Let $h\in[g_1\cup g_2]$. If $h\cap g_1$ is nonempty-valued, then it is usco by Lemma~\ref{usco-stability}. Then it follows from the fact that $[g_1]=\{f\}$ that $f\subset h\cap g_1$. By the minimality of $h$ we get $f=h$.

Therefore, if $h\ne f$, there is $x\in X$ such that $h(x)\cap g_1(x)=\emptyset$. Using the Hausdorff property of $Y$ we get an open set $V\subset X$ containing $x$ such that $h(V)\cap g_1(V)=\emptyset$. Define  a map $\tilde h:X\to Y$ by the formula
$$\tilde h(x)=\begin{cases} g_2(x), & x\in X\setminus V, \\ h(x), & x\in V.\end{cases}$$
Then $\tilde h\subset g_2$ and by Lemma~\ref{usco-modif} it is usco.
Hence $f\subset \tilde h$, a contradiction.
\end{proof}

\begin{lemma}\label{q-uniq}
Let $X$ be a Baire topological space, $Y$ a metrizable space and $g:X\to Y$ a quasiminimal usco mapping. Then $g(x)$ is a singleton for all $x$ in a residual subset of $X$.
\end{lemma}

\begin{proof} The proof is done by a minor modification of the proof of \cite[Proposition 3.1.4]{fabian}.
Denote by $f$ the unique element of $[g]$. Let $\rho$ be a metric
generating the topology of $Y$. For $n\in\mathbb{N}$ set
$$G_n=\{x\in X: \exists U\subset X\mbox{ open}, x\in U: \rho\mbox{-}\diam g(U)\le\frac1n\}.$$
Then each $G_n$ is clearly open and for any $x\in\bigcap_{n\in\mathbb{N}} G_n$ the image $g(x)$ is a singleton. Hence it is sufficient to show that each $G_n$ is dense.

To see it let $V\subset X$ be an arbitrary nonempty open set.
Choose any point $y\in f(V)$. By Lemma~\ref{minimal-prop} there is
a nonempty open subset $W\subset V$ with $f(W)\subset
B_\rho(y,\frac1{2n})$. We claim that there is some $x\in W$ with
$g(x)\subset B_\rho(y,\frac1{2n})$. Indeed, suppose it is not the
case. Then the mapping
$$h(x)=\begin{cases} g(x), & x\in X\setminus W, \\ g(x)\setminus B_\rho(y,\frac1{2n}), & x\in W,
\end{cases}$$
is usco by Lemma~\ref{usco-modif}. As $[g]=\{f\}$, necessarily
$f\subset h$, which is a contradiction as $f(x)\cap h(x)=\emptyset$
for $x\in W$.

Finally, as $g$ is usco, there is an open neighborhood $U$ of $x$
with $g(U)\subset B_\rho(y,\frac1{2n})$ and hence $x\in V\cap
G_n$.
\end{proof}

\begin{remark} The assertion of Lemma~\ref{q-uniq} remains true if one assumes, instead of metrizability of $Y$, that
the space $Y$ is fragmented by some metric. Let us recall the
definition of fragmentability. If $(Y,\mathcal{T})$ is a topological
space and $\rho$ is a metric on the set $Y$, we say that $Y$ is {\it
fragmented by $\rho$\/} if each nonempty subset of $Y$ has nonempty
relatively $\mathcal{T}$-open subsets of arbitrarily small
$\rho$-diameter. The proof of the more general statement is almost
the same as the proof of Lemma~\ref{q-uniq}, one should only modify
the proof of \cite[Theorem 5.1.11]{fabian} instead of that of
\cite[Proposition 3.1.4]{fabian}. However, the following question
seems to be open.

\begin{question} Is any quasiminimal usco from a Baire space into a
Stegall space singlevalued at points of a residual set?
\end{question}

Recall that a space $Y$ is {\it Stegall} if any minimal usco from a
Baire space into $Y$ is singlevalued at points of a residual set
(see e.g. \cite[Chapter 3]{fabian}).
\end{remark}

\section{Convergence structure on $\mathcal{M}(X,Y)$}\label{sec-def}

In this section we will define convergence of filters on $\mathcal{M}(X,Y)$ and show that it defines a convergence structure.

\begin{definition}\label{defFilterConv}
A filter $\mathcal{F}$ on $\mathcal{M}(X,Y)$ {\it converges to
$f\in\mathcal{M}(X,Y)$} and we write $\mathcal{F}\rightarrow f$ if
$$\{f\}=\bigcap\{ [g] : g\mbox{ is a usco map}, [g]\in\mathcal{F} \}.$$
\end{definition}

\begin{remarks}\label{remDefinition}\quad
\begin{itemize}
\item[1.] If $\mathcal{F}\rightarrow f$, then there is at least
one usco map $g:X\to Y$ with $[g]\in \mathcal{F}$. If the filter
$\mathcal{F}$ is such that there exists an usco map $g:X\to Y$ with
$[g]\in \mathcal{F}$ it is called {\it usco-bounded}. Since this is
the only concept of boundedness considered in the paper we call the
usco-bounded filters simply bounded.

\item[2.] If $g_1$ and $g_2$ are two usco maps such that both
$[g_1]$ and $[g_2]$ belong to  a filter $\mathcal{F}$, then
$[g_1]\cap[g_2]\in\mathcal{F}$ as well. Thus $g_1\cap g_2$ is
nonempty-valued and therefore it is a usco map by
Lemma~\ref{usco-stability}. Further, clearly $[g_1\cap
g_2]=[g_1]\cap[g_2]$. It follows that the family
\begin{equation}\label{familyg}
\{ [g]: g\mbox{ is a usco map}, [g]\in\mathcal{F} \}
\end{equation}
is closed to finite intersections and hence it is a filter base
provided it is nonempty. For any bounded filter $\mathcal{F}$  on
$\mathcal{M}(X,Y)$ we denote by $\mathcal{G}_\mathcal{F}$ the filter
which is generated by the family (\ref{familyg}). Obviously,
$\mathcal{G}_\mathcal{F}$ is coarser then $\mathcal{F}$ and we have
from the definition that
\begin{equation}\label{GFconv}
\mathcal{F}\rightarrow f\ \Longleftrightarrow\
\mathcal{G}_{\mathcal{F}}\rightarrow f
\end{equation}

\item[3.] Let $\mathcal{F}$ be a bounded filter. Set
   $$ g_{\mathcal{F}}=\bigcap\{ g : g\mbox{ is a usco map}, [g]\in\mathcal{F} \}.$$
   Then $g_{\mathcal{F}}$ is a usco map (by the previous remark and Lemma~\ref{usco-stability}) and we have
\[
\bigcap_{F\in\mathcal{G}_\mathcal{F}} F=[g_\mathcal{F}].
\]
Further, $\mathcal{F}\to f$ if and only if $[g_{\mathcal{F}}]=\{f\}$.

\item[4.] It is obvious from the definition that one filter cannot
converge to more than one element of $\mathcal{M}(X,Y)$.

\item[5.] If the domain space $X$ is a singleton, then
$\mathcal{M}(X,Y)$ can be canonically identified with $Y$. Then a
filter $\mathcal{F}$ on $Y$ converges to $y\in Y$ in
$\mathcal{M}(X,Y)$ if and only if it contains a compact subset of
$Y$ and converges to $y$ in the topology of $Y$.
\end{itemize}
\end{remarks}

According to Definition \ref{defFilterConv} with every point
$f\in\mathcal{M}(X,Y)$ we associate a set of filters $\lambda(f)$
which converge to $f$. The mapping $\lambda$ from
$\mathcal{M}(X,Y)$ into the power set of the set of filters on
$\mathcal{M}(X,Y)$ is called a convergence structure and
$(\mathcal{M}(X,Y),\lambda)$ is called a convergence space if the
following conditions are satisfied for all $f\in\mathcal{M}(X,Y)$,
see \cite{Beattie}:
\begin{eqnarray}
\bullet&&\langle f\rangle\in\lambda(f),\ \mbox{ where $\langle
f\rangle$ denotes the filter
generated by $\{\{f\}\}$}.\label{ConStrCond1}\\
\bullet&&\mbox{If $\mathcal{F}_1,\mathcal{F}_2\in\lambda(f)$ then
$\mathcal{F}_1\cap\mathcal{F}_2\in\lambda(f)$.}\label{ConStrCond2} \\
\bullet&&\mbox{If $\mathcal{F}_1\in\lambda(f)$ then
$\mathcal{F}_2\in\lambda(f)$ for all filters $\mathcal{F}_2$ on }
\mathcal{M}(X,Y)\nonumber\\
&&\mbox{which are finer than }\mathcal{F}_1.\label{ConStrCond3}
\end{eqnarray}

\begin{theorem}\label{tConvSpace}
The mapping $\lambda$ is a convergence structure on
$\mathcal{M}(X,Y)$.
\end{theorem}

\begin{proof} We need to show
that for every $f\in\mathcal{M}(X,Y)$ conditions
(\ref{ConStrCond1})--(\ref{ConStrCond3}) are satisfied. Conditions
(\ref{ConStrCond1}) and (\ref{ConStrCond3}) follow immediately from
Definition~\ref{defFilterConv}. We will show that condition
(\ref{ConStrCond2}) also holds. Let
$\mathcal{F}_1,\mathcal{F}_2\in\lambda(f)$.  We define the following
set of usco mappings:
\[
\Phi=\{g^{(1)}\cup g^{(2)}:g^{(1)},g^{(2)}\mbox{ are usco maps},
[g^{(1)}]\in\mathcal{F}_1,[g^{(2)}]\in\mathcal{F}_2\}.
\]
By Lemma~\ref{usco-stability} the family $\Phi$ consists of usco
maps. As
$$ [g^{(1)}\cup g^{(2)}]\supset [g^{(1)}]\cup[g^{(2)}],$$
we get $\{[h]:h\in\Phi\}\subset\mathcal{F}_1\cap\mathcal{F}_2$. Set
$g=\bigcap\Phi$. It is easy to check that
$$g=g_{\mathcal{F}_1}\cup g_{\mathcal{F}_2}.$$
As $[g_{\mathcal{F}_1}]=[g_{\mathcal{F}_2}]=\{f\}$, it follows by
Lemma~\ref{union} that $[g]=\{f\}$ as well. From the inclusion
$[g_{\mathcal{F}_1\cap\mathcal{F}_2}]\subseteq[g]$ it follows that
$[g_{\mathcal{F}_1\cap\mathcal{F}_2}]=\{f\}$. Hence
$\mathcal{F}_1\cap\mathcal{F}_2\to f$, which completes the proof.
\end{proof}

Let $(f_\nu)_{\nu\in I}$ be a net in $\mathcal{M}(X,Y)$ indexed by a directed set $I$.
 Following the general theory of convergence spaces the net
$(f_\nu)_{\nu\in I}$ converges to $f\in\mathcal{M}(X,Y)$ if
 the filter generated by
$\{\{f_\nu:\nu\geq \nu_0\}:\nu_0\in I\}$, converges to $f$ in the
convergence structure $\lambda$.
 In particular, a sequence
$(f_n)_{n\in\mathbb{N}}$ converges to $f\in\mathcal{M}(X,Y)$ if
its Fr\'echet filter, that is, the filter generated by
$\{\{f_m:m\geq n\}:n\in\mathbb{N}\}$, converges to $f$ in the
convergence structure $\lambda$. The following theorem gives an
alternative characterization of the convergent nets and sequences.

\begin{theorem}\label{tConvSeq}\quad
\begin{itemize}
    \item[a)] A net $(f_\nu)_{\nu\in I}$ converges to $f\in\mathcal{M}(X,Y)$ if and only if there is some $\nu_0\in I$ and usco mappings $g_\nu$, $\nu\ge\nu_0$ such that

\begin{itemize}
    \item[(i)] $f_\nu\subset g_\nu$ for $\nu\ge\nu_0$;
    \item[(ii)] $g_\nu\subset g_{\nu'}$ for $\nu\ge\nu'\ge\nu_0$;
    \item[(iii)] $f$ is the unique minimal usco contained in $\bigcap\limits_{\nu\ge\nu_0} g_\nu$.
\end{itemize}
\item[b)] A sequence $(f_n)_{n\in\mathbb{N}}$ converges to
$f\in\mathcal{M}(X,Y)$ if and only if there exists a sequence of
usco maps $(g_n)_{n\in\mathbb{N}}$ such that
\begin{itemize}
\item[(i)] $f_n\subset g_n$ for $n\in\mathbb{N}$; \item[(ii)]
$g_m\subset g_n$ for each $m\ge n$; \item[(iii)] $f$ is the unique
minimal usco contained in $\bigcap\limits_{n\in\mathbb{N}}g_n$.
\end{itemize}
\end{itemize}
\end{theorem}

\begin{proof}
a) Let $(f_\nu)_{\nu\in I}$ converge to $f\in\mathcal{M}(X,Y)$.
Denote by $\mathcal{F}$ the filter generated by $\{\{f_\nu:\nu\geq
\nu_0\}:\nu_0\in I\}$. For $\nu\in I$ set
$$g_\nu=\bigcap\{g : g\mbox{ is a usco map}, f_{\nu'}\subset g\mbox{ for }\nu'\ge\nu\}.$$
As $\mathcal{F}\to f$, there is a usco map $g$ such that
$[g]\in\mathcal{F}$. Then there is some $\nu_0$ such that
$f_\nu\subset g$ for $\nu\ge\nu_0$. Therefore $g_{\nu_0}$ is a
well defined usco map (by Lemma~\ref{usco-stability}). Hence
$g_\nu$ is a well defined usco for each $\nu\ge\nu_0$. The
conditions (i) and (ii) are satisfied by the definition. To see
that the condition (iii) is satisfied too, it suffices to observe
that for any usco map $g$ with $[g]\in\mathcal{F}$ there is
$\nu\ge\nu_0$ with $g_\nu\subset g$.

Assume now that there exists a net of usco maps
$(g_\nu)_{n\ge\nu_0}$ satisfying conditions (i), (ii) and (iii).
It follows from (i) and (ii) that $[g_\nu]\in\mathcal{F}$ for each
$n\ge n_0$. Thus, due to (iii), $\mathcal{F}\to f$.

b) Suppose that $f_n\to f$. It follows from (1) that there is
$n_0\in\mathbb{N}$ and usco maps $g_n$, $n\ge n_0$, satisfying a),
b) and c). For $n<n_0$ we take
\[
g_n=f_n\cup ... \cup f_{n_0-1} \cup g_{n_0}.
\]
Then usco maps $g_n$, $n\in\mathbb{N}$, fulfil the conditions (i),
(ii) and (iii).

The inverse implication follows from that in a).
\end{proof}

\begin{remark}
The preceding theorem indicates a relation of the convergence on
$\mathcal{M}(X,Y)$ to  the order convergence on a lattice. We will
examine the relationship in more detail in Section \ref{sec-order}.
\end{remark}

The preceeding theorem enables us to show that the convergence on $\mathcal{M}(X,Y)$ is not in general generated by a topology.

\begin{example}\label{nontop} If $X=Y=[0,1]$, then the convergence in $\mathcal{M}(X,Y)$ is not generated by any topology.
\end{example}

\begin{proof} Let $q_n$, $n\in\mathbb{N}$ be an enumeration of rational numbers from $[0,1]$. We define continuous functions $f_n:[0,1]\to[0,1]$ by the formula
$$f_n(x)=\begin{cases} 1-n|x-q_n|, & x\in (q_n-\frac1n,q_n+\frac1n)\cap[0,1], \\0, & \mbox{otherwise}.\end{cases}$$
Then the sequence $f_n$ does not converge in $\mathcal{M}(X,Y)$. Indeed, if $\mathcal{F}$ is the Fr\'echet filter of this sequence, it is easy to check that
$$g_{\mathcal{F}}(x)=[0,1], \qquad x\in[0,1],$$
which is obviously not quasiminimal.

On the other hand, if a subsequence $q_{n_k}$ converges to some $q\in[0,1]$, the sequence $f_{n_k}$ converges to $0$. Indeed, if we set
$$g_k(x)=\begin{cases} \{ f_{n_l}(x): l\ge k\}  , & x\in [0,1]\setminus\{q\}, \\
[0,1], & x=q,\end{cases}$$
we get a decreasing sequence of usco maps such that $g_k\supset
f_{n_k}$ for each $k$. Further,
$$\bigcap_{k\in\mathbb{N}} g_k(x)=\begin{cases} [0,1], & x=q, \\ \{0\}, & x\in[0,1]\setminus \{q\}. \end{cases}$$
This usco map is clearly quasiminimal and the only minimal usco contained in it is the constant zero function. Thus $f_{n_k}$ converges to $0$ by Theorem~\ref{tConvSeq}.

As each subsequence of $q_n$ has a further convergent subsequence, we get that each subsequence of $f_n$ has a further subsequence converging to $0$. If the convergence were a topological one,
it would imply that $f_n$ converge to $0$ as well. But it is not the case by the first paragraph, hence the convergence is not a topological one.
\end{proof}

\begin{theorem}
The convergence given in Definition \ref{defFilterConv} is stable
with respect to restrictions to open sets, that is, if a filter
$\mathcal{F}$ on $\mathcal{M}(X,Y)$ converges to
$f\in\mathcal{M}(X,Y)$ then for any open subset $D$ of $X$ the
filter $\mathcal{F}|_D$ generated by the restrictions $\{F|_D:F\in
\mathcal{F}\}$, where $F|_D=\{h|_D:h\in F\}\subset
\mathcal{M}(D,Y)$, converges on $D$ to the restriction of the limit
$f|_D$.
\end{theorem}

\begin{proof}
The restriction of a minimal usco to an open set is also a minimal
usco, \cite[Lemma 2]{stenfra}. Therefore $\mathcal{F}|_D$ is a
filter on $\mathcal{M}(D,Y)$ and $f|_D\in\mathcal{M}(D,Y)$.

 Set
 $$\Phi = \{ g|_D : g\mbox{ is a usco map}, [g]\in\mathcal{F} \}.$$
 Then $\Phi$ is a family of usco maps from $D$ to $Y$. Further, $$[\Phi]=\{[h]: h\in \Phi\}\subset\mathcal{F}|_D.$$
 Indeed, we have
 $$[g]|_D=\{h|_D: h\in[g]\}\subset [g|_D]$$
 for any usco map $g$.
 Obviously $f|_D\in\bigcap[\Phi]$. It remains to show that $\bigcap[\Phi]$ has no more elements.
 Let $h\in\bigcap[\Phi]$ be any element. Set
 $$\psi(x)=\begin{cases} h(x), & x\in D, \\ g_{\mathcal F}(x), & x\in X\setminus D.\end{cases}$$
 Then $\psi$ is usco by Lemma~\ref{usco-modif}. Further, $\psi\subset g_{\mathcal{F}}$, and hence $f\subset\psi$. It follows that $f|_D\subset h$, hence $f|_D=h$.
 This completes the proof.
 \end{proof}

\section{Relationship to pointwise convergence}\label{sec-pointwise}

In this section we give some relations of the convergence in $\mathcal{M}(X,Y)$ to the pointwise convergence.

\begin{theorem}\label{pointwise} Let $X$ be a Baire space, $f_n$ be a bounded sequence in $\mathcal{M}(X,Y)$ and $\varphi:X\to Y$ be a quasiminimal usco with $[\varphi]=\{f\}$. Suppose that for each $x\in X$ the sequence of compact sets $f_n(x)$ {\it cumulates at} $\varphi(x)$, i.e. for each open set $U\subset Y$ containing $\varphi(x)$ there is some $n_0\in\mathbb{N}$ such that $f_n(x)\subset U$ for $n\ge n_0$. Then the sequence $f_n$ converges to $f$ in $\mathcal{M}(X,Y)$.
\end{theorem}

 \begin{proof} Let $g_n$ be the intersection of all usco maps containing $f_k$ for $k\ge n$ and let $g$ be the intersection of all $g_n$'s. As the sequence $f_n$ is bounded, $g_n$'s and $g$ are well-defined usco maps. By Theorem~\ref{tConvSeq} it suffices to prove that $[g]=\{f\}$.

Choose $h\in [g]$ arbitrary. Suppose that $h\ne f$. Then there is some $x_0\in X$ such that $h(x_0)\cap \varphi(x_0)=\emptyset$. Indeed, otherwise $h\cap \varphi$ would be a usco map contained in $h$  and hence we would have $h \cap\varphi=h$, i.e. $h\subset\varphi$. But then necessarily $h=f$.

As $Y$ is Hausdorff, there are disjoint open sets $V_1\supset \varphi(x_0)$ and $V_2\supset h(x_0)$.
Further, there is an open neighborhood $U$ of $x_0$ such that $\varphi(U)\subset V_1$ and $h(U)\subset V_2$.
For each $n\in\mathbb{N}$ set
$$ F_n=\{ x\in X: (\forall k\ge n)(f_k(x)\setminus V_2\ne\emptyset) \}.$$
Then each $F_n$ is a closed subset of $X$. Moreover, the sets $F_n$
cover $U$. Indeed, if $x\in U$, then $\varphi(x)\subset V_1$. Hence
there is some $n$ such that for each $k\ge n$ we have $f_k(x)\subset
V_1$. Thus $x\in F_n$. As $X$ is a Baire space, $U$ is non-meager
and hence there is some $n\in\mathbb{N}$ such that $F_{n}\cap U$ has
nonempty interior. It means that there is a nonempty open set
$W\subset U$ such that $f_k(x)\setminus V_2\ne\emptyset$ whenever
$k\ge n$ and $x\in W$. Hence, by the minimality of $f_k$ we get
(using Lemma~\ref{minimal-prop}) that $f_k(W)\cap V_2=\emptyset$ for
each $k\ge n$. Therefore, if we define
$$\widetilde g_n(x)=\begin{cases} g_n(x)\setminus V_2, & x\in W, \\ g_n(x), & x\in X\setminus W,
\end{cases}$$
we get a usco map containing $f_k$ for $k\ge n$. Thus $\widetilde g_n= g_n$. As $g\subset g_n$, we get that $g(x)\cap V_2=\emptyset$ for each $x\in W$, which is a contradiction with the assumption $h\in[g]$. This completes the proof.
\end{proof}

As a corollary we get the following assertions on sequences of continuous functions.

\begin{corollary}\label{corseq} Let $X$ be a Baire space, $f_n$ be a sequence of continuous functions bounded in $\mathcal{M}(X,Y)$.
\begin{itemize}
    \item[(i)] If $f_n$ pointwise converges to a continuous function $f$, then $f_n$ converges to $f$ in $\mathcal{M}(X,Y)$.
    \item[(ii)] If $f\in\mathcal{M}(X,Y)$ is such that the sequence $f_n(x)$ converges to an element of $f(x)$ for each $x\in X$, then $f_n$ converges to $f$ in $\mathcal{M}(X,Y)$.
\end{itemize}
\end{corollary}

The following example shows that all assumptions in Theorem~\ref{pointwise} are needed. Namely, one can drop neither the assumption that the sequence is bounded (even if $X$ is compact),
nor the assumption that $X$ is a Baire space (even if $Y$ is compact and hence all filters are bounded). Further, Theorem~\ref{pointwise} is not true for nets, even if both $X$ and $Y$ are
compact.

\begin{example}\label{pw-exa}\quad
\begin{itemize}
    \item[1.] There is a sequence of continuous functions $f_n:[0,1]\to\mathbb{R}$ pointwise converging to $0$ which is unbounded in $\mathcal{M}([0,1],\mathbb{R})$.
    \item[2.] There is a sequence of continuous functions $f_n:\mathbb{Q}\to[0,1]$ pointwise converging to $0$ which is not convergent in $\mathcal{M}(\mathbb{Q},[0,1])$.
    \item[3.] There is a net of continuous functions $f_\nu:[0,1]\to[0,1]$ pointwise converging to $0$ which is not convergent in $\mathcal{M}([0,1],[0,1])$.
\end{itemize}
\end{example}

\begin{proof} 1. It is sufficient to take the sequence
$$f_n(x)= \begin{cases} n(1-2n|x-\frac1{2n}|), & x\in [0,\frac1n], \\
0 & \mbox{otherwise.}\end{cases}$$

2. Let $q_n$, $n\in\mathbb{N}$, be an enumeration of $\mathbb{Q}$.
For each $n\in\mathbb{N}$ choose a continuous function
$f_n:\mathbb{Q}\to[0,1]$ such that
$$f_n(x)=\begin{cases} 0 & x\in\{q_1,\dots,q_n\}, \\
1 & x\in\mathbb{Q}\setminus\bigcup_{k=1}^n (q_k-\frac1{n^2},q_k+\frac1{n^2}).\end{cases}$$
Then $f_n$ pointwise converge to $0$. Further, fix $n\in\mathbb{N}$ and a usco map $g:\mathbb{Q}\to[0,1]$ containing $f_m$ for $m\ge n$. Then $0\in g(x)$ for each $x\in\mathbb{Q}$. Moreover, $1\in g(x)$ for each
$$x\in \bigcup_{m\ge n}\left(\mathbb{Q}\setminus \bigcup_{k=1}^m \left(q_k-\frac1{m^2},q_k+\frac1{m^2}\right)\right)
=\mathbb{Q}\setminus \bigcap_{m\ge n}\bigcup_{k=1}^m
\left(q_k-\frac1{m^2},q_k+\frac1{m^2}\right).
$$
As $\bigcap_{m\ge n}\bigcup_{k=1}^m [q_k-\frac1{m^2},q_k+\frac1{m^2}]$ is a closed subset of $\mathbb{R}$ of Lebesgue measure zero, its complement is an open dense set. Therefore $1\in g(x)$ for all $x$ from a dense subset of $\mathbb{Q}$, and thus for all $x\in\mathbb{Q}$.

It follows that $\{0,1\}\subset g(x)$ for each $x\in\mathbb{Q}$ and hence the sequence $f_n$ is not convergent in $\mathcal{M}(\mathbb{Q},[0,1])$ by Theorem~\ref{tConvSeq}.

3. For any nonempty finite set $A\subset[0,1]$ choose a continuous
function $f_A:[0,1]\to[0,1]$ such that
$$f_A(x)=\begin{cases} 0 & x\in A, \\
1 & x\in[0,1]\setminus\bigcup_{a\in A}
\left(a-\frac1{|A|^2},a+\frac1{|A|^2}\right),\end{cases}$$ where
$|A|$ denotes the cardinality of $A$. If we consider finite subsets
of $[0,1]$ ordered by inclusion, the net $f_A$ pointwise converges
to $0$. Moreover, the net $f_A$ is not convergent in
$\mathcal{M}([0,1],[0,1])$. Indeed, let $B\subset [0,1]$ be a
nonempty finite set and $g:[0,1]\to[0,1]$ be a usco map containing
$f_A$ for each $A\supset B$. Then $0\in g(x)$ for each $x\in [0,1]$.
Further, $1\in g(x)$ for each
$$x\in\bigcup_{A\supset B} \left([0,1]\setminus\bigcup_{a\in A} \left(a-\frac1{|A|^2},a+\frac1{|A|^2}\right)\right) =
[0,1]\setminus\bigcap_{A\supset B}\bigcup_{a\in A}
\left(a-\frac1{|A|^2},a+\frac1{|A|^2}\right)$$
Again, $\bigcap_{A\supset B}\bigcup_{a\in A} (a-\frac1{|A|^2},a+\frac1{|A|^2})$ is a set of Lebesgue measure zero, hence the complement is dense in $[0,1]$. Therefore, $1\in g(x)$ for all $x$ from a dense subset of $[0,1]$ and thus for all $x\in[0,1]$. It follows that $\{0,1\}\subset g(x)$ for all $x\in[0,1]$ and so the net $f_A$ does not converge by Theorem~\ref{tConvSeq}.
\end{proof}

In the next theorem we give a result in the converse direction.
Let $F\subset\mathcal{M}(X,Y)$. For a given $x\in X$, $F(x)$ is
the set
\[
F(x)=\bigcup\{f(x):f\in F\}
\]
Let $\mathcal{F}$ be a filter on $\mathcal{M}(X,Y)$. Then for
every $x\in X$
\[
\{F(x):F\in \mathcal{F}\}
\]
is a filter base on $Y$. Denote by $\mathcal{F}(x)$ the filter it
generates. We consider the following question. Given that a filter
$\mathcal{F}$ converges to $f\in\mathcal{M}(X,Y)$, is the filter
$\mathcal{F}(x)$ convergent in the topology of $Y$ for some $x\in
X$? The next theorem deals with this question in the case when $X$
is a Baire space and $Y$ is a metric space.
%Let us recall
%that in this case the minimal usco maps from $X$ to $Y$ are single
%valued and continuous on a residual subset of the domain.

\begin{theorem}\label{tPointConv}
Let $X$ be a Baire space and $Y$ be a metric space with metric
$\rho$. If a filter $\mathcal{F}$ on $\mathcal{M}(X,Y)$ converges
to $f\in\mathcal{M}(X,Y)$ then there exists a residual subset $D$
of $X$ such that for every $x\in D$
\begin{itemize}
\item[(i)] $f(x)$ is a singleton, that is, $f(x)\in Y$;

\item[(ii)] $\mathcal{F}(x)$ converges to $f(x)$ with respect to
the metric $\rho$.
\end{itemize}
\end{theorem}

\begin{proof}
As $\mathcal{F}\to f$, the usco map $g_{\mathcal F}$ is quasiminimal and $[g_{\mathcal F}]=\{f\}$.
By Lemma~\ref{q-uniq} there is a residual set $D\subset X$ such that $g_{\mathcal{F}}(x)=f(x)$ is a singleton for each $x\in D$.

Fix $x\in D$. By the previous paragraph we have
$$\{f(x)\} =\bigcap\{ g(x): g\mbox{ is a usco map}, [g]\in\mathcal{F} \}.$$
If $g$ is a usco map such that $[g]\in\mathcal{F}$, then
$$g(x)\supset \bigcup\{h(x):h\in[g]\}=[g](x)\in\mathcal{F}(x).$$
As these $g(x)$'s are compact subsets of $Y$ belonging to $\mathcal{F}(x)$
and their intersection is just $\{f(x)\}$, we get $\mathcal{F}(x)\to f(x)$.
This completes the proof.
\end{proof}

\section{Convergence in $\mathcal{M}(X,\mathbb{R})$ and the order
convergence}\label{sec-order}

In this section we show that the convergence in $\mathcal{M}(X,\mathbb{R})$ is equivalent to the order convergence with respect to the natural partial order on $\mathcal{M}(X,\mathbb{R})$.
Before giving the definitions and stating the equivalence we show a natural correspondence of the space $\mathcal{M}(X,\mathbb{R})$ and the space $\mathbb{H}(X,\mathbb{R})$ of Hausdorff continuous functions (see \cite{Sendov,Anguelov}).

We start by the following obvious lemma.

\begin{lemma}\label{convex}\quad
\begin{itemize}
    \item Let $g:X\to\mathbb{R}$ be a usco map. Then the map $x\mapsto \max g(x)$ is upper semicontinuous on $X$ and the map $x\mapsto \min g(x)$ is lower semicontinuous on $X$.
    \item Let $f_1:X\to\mathbb{R}$ be a lower semicontinuous function and $f_2:X\to\mathbb{R}$ be an upper semicontinuous function such that $f_1\le f_2$ on $X$. Then the set-valued map $x\mapsto[f_1(x),f_2(x)]$ is usco.
\end{itemize}
\end{lemma}

%We consider next a natural connection which exists between the
%convergence according to Definition \ref{defFilterConv} and the
%order convergence when $Y=\mathbb{R}$. In this case a partial
%order is induced on $\mathcal{M}(X,\mathbb{R})$ by the linear
%order on $\mathbb{R}$.

Let $f\in\mathcal{M}(X,\mathbb{R})$. We define the following two real functions on $X$:
\begin{eqnarray*}
\underline{f}(x)&=&\min f(x)\\
\overline{f}(x)&=&\max f(x)
\end{eqnarray*}
By the previous lemma $\underline{f}$ is lower semicontinuous and $\overline{f}$ is upper semicontinuous. Further, we define a map $f_C:X\to\mathbb{R}$ by
$$ f_C(x)=[\underline{f}(x),\overline{f}(x)].$$
By Lemma~\ref{convex} it is a usco map. Moreover, we have the following.

\begin{lemma}\label{f-C}\quad
\begin{itemize}
    \item[(i)] If $f,g\in\mathcal{M}(X,\mathbb{R})$ are distinct, then $f_C(x)\cap g_C(x)=\emptyset$ for some $x\in X$.
    \item[(ii)] For each $f\in\mathcal{M}(X,\mathbb{R})$ we have $[f_C]=\{f\}$.
    \item[(iii)] For each $f\in\mathcal{M}(X,\mathbb{R})$ the usco map $f_C$ is minimal within the convex valued (i.e., interval-valued) usco maps.
\end{itemize}
\end{lemma}

\begin{proof}
(i) Let $f$ and $g$ be distinct elements of
$\mathcal{M}(X,\mathbb{R})$. Then there is some $x_0\in X$ with
$f(x_0)\cap g(x_0)=\emptyset$ (otherwise $f\cap g$ would be a usco
contained both in $f$ and $g$ and hence we would have $f=g=f\cap
g$). Let $a=\max g(x_0)$ and $b=\max f(x_0)$. Then $a\ne b$, we can
suppose without loss of generality that $a<b$. Choose some $c\in
(a,b)$. As $g$ is usco, there is an open neighborhood $U$ of $x_0$
such that $g(U)\subset (-\infty,c)$. We have $f(U)\cap
(c,+\infty)\ne\emptyset$, and hence by Lemma~\ref{minimal-prop}
there is a nonempty open set $V\subset U$ with $f(V)\subset
(c,+\infty)$. It follows that $f_C(x)\cap g_C(x)=\emptyset$ for any
$x\in V$.

(ii) This follows immediately from (i).

(iii) Let $g\subset f_C$ be an interval-valued usco. By (ii) we have $f\in[g]$. Hence it follows from the definition of $f_C$ that $f_C\subset g$.
\end{proof}

It is easy to check that the minimal interval-valued usco maps are
exactly the Hausdorff continuous functions in the sense of
\cite{Sendov}. Hence, due to the previous lemma the mapping
$f\mapsto f_C$ is a bijection of $\mathcal{M}(X,\mathbb{R})$ onto
$\mathbb{H}(X,\mathbb{R})$. On the set $\mathbb{H}(X,\mathbb{R})$
there is a natural partial order (see \cite{Anguelov}). We define
a partial order on $\mathcal{M}(X,\mathbb{R})$ using the
correspondence $f\mapsto f_C$:

 For $f,g\in\mathcal{M}(X,\mathbb{R})$ we have
\begin{equation}\label{order}
f\leq g\ \Longleftrightarrow \
\underline{f}(x)\leq\underline{g}(x),\
\overline{f}(x)\leq\overline{g}(x),\ x\in X.
\end{equation}

Using the minimality of $f$ and $g$ it is easy to see that either
one of the inequalities on the right hand side above will suffice,
that is, we have
\begin{equation}\label{order1}
f\leq g\Longleftrightarrow \underline{f}(x)\leq\underline{g}(x),\
x\in X\Longleftrightarrow \overline{f}(x)\leq\overline{g}(x),\
x\in X.
\end{equation}
Indeed, let $\underline{f}\le\underline{g}$ on $X$. The function $h(x)=\min\{\overline{f}(x),\overline{g}(x)\}$ is upper semicontinuous and clearly
$$\underline{f}\le h \le \overline{f}$$
on $X$. Hence the map $x\mapsto [\underline{f}(x),h(x)]$ is an interval-valued usco
(see Lemma~\ref{convex}) contained in $f_C$. Then it is equal to $f_C$ by Lemma~\ref{f-C}.
Hence $h=\overline{f}$, which means $\overline{f}\le\overline{g}$ on $X$. This proves one
implication, the inverse one can be proved in the same way.

Since the mapping $f\mapsto f_C$ is an order isomorphism of from
$\mathcal{M}(X,\mathbb{R})$ onto $\mathbb{H}(X,\mathbb{R})$ the set
$\mathcal{M}(X,\mathbb{R})$ has the same order properties as
$\mathbb{H}(X,\mathbb{R})$. For example, since
$\mathbb{H}(X,\mathbb{R})$ is Dedekind order complete, see
\cite{Anguelov}, $\mathcal{M}(X,\mathbb{R})$ is also Dedekind order
complete. In particular this implies that
$\mathcal{M}(X,\mathbb{R})$ is a lattice.

The following theorem shows an essential similarity between the
the functions in $\mathcal{M}(X,\mathbb{R})$ and the usual
continuous real valued functions on $X$. It follows from the
respective statement for Hausdorff continuous functions, see
\cite[Theorem 4]{Anguelov}.

\begin{theorem}\label{tdense}
Let $f,g\in\mathcal{M}(X,\mathbb{R})$ and let $D$ be a dense
subset of $X$. Then
\[
f|_D\leq g|_D \Longrightarrow f\leq g
\]
\end{theorem}

Next we will establish a link between the order convergence on
$\mathcal{M}(X,\mathbb{R})$ with respect to the order
(\ref{order}) and the convergence structure $\lambda$. Let us
recall the definition for order convergence of filters. For a
filter $\mathcal{F}$ on $\mathcal{M}(X,\mathbb{R})$ we consider
the set of lower bounds
\[
\mathcal{F}^-=\{\phi\in\mathcal{M}(X,\mathbb{R}):\exists
F\in\mathcal{F}:\phi\leq h \mbox{ for all } h\in F \}
\]
and the set of upper bounds
\[
\mathcal{F}^+=\{\psi\in\mathcal{M}(X,\mathbb{R}):\exists
F\in\mathcal{F}:\psi\geq h \mbox{ for all } h\in F \}.
\]
We say that the filter $\mathcal{F}$ order converges to
$f\in\mathcal{M}(X,\mathbb{R})$ if $\mathcal{F}^-\neq\emptyset$,
$\mathcal{F}^+\neq\emptyset$ and
\begin{equation}\label{f=sup=inf}
f=\sup \mathcal{F}^-=\inf \mathcal{F}^+.
\end{equation}

\begin{remark} Let us notice that $\phi\in\mathcal{F}^-$ and
$\psi\in\mathcal{F}^+$ if and only if the order interval
$[\phi,\psi]$ belongs to $\mathcal{F}$. Hence a filter $\mathcal{F}$
order converges to $f$ if and only if
$$\{f\}=\bigcap\{ [\psi,\phi] :
\psi,\phi\in\mathcal{M}(X,\mathbb{R}),[\psi,\phi]\in\mathcal{F} \}
$$
\end{remark}

\begin{theorem}\label{tOrderConv}
A filter $\mathcal{F}$ on $\mathcal{M}(X,\mathbb{R})$ order
converges to $f\in\mathcal{M}(X,\mathbb{R})$ iff
$\mathcal{F}\in\lambda(f)$.
\end{theorem}
\begin{proof}
Let $\mathcal{F}$ order converge to $f$. For arbitrary $\phi\in
\mathcal{F}^-$ and $\psi\in \mathcal{F}^+$ we have $\phi\leq\psi$.
Hence the usco map
$h_{\psi,\phi}:x\mapsto[\underline{\phi}(x),\overline{\psi}(x)]$ is
well defined on $X$ and $[h_{\psi,\phi}]\in\mathcal{F}$. Due to
(\ref{f=sup=inf}) we have that $f$ is the only minimal usco
contained in the map
\[
\varphi=\bigcap\{h_{\phi,\psi}:\phi\in \mathcal{F}^-,\psi\in
\mathcal{F}^+\}  :x\mapsto\ \bigcap_{\phi\in \mathcal{F}^-,\,\psi\in
\mathcal{F}^+}[\underline{\phi}(x),\overline{\psi}(x)],\ x\in X.
\]
Since $g_\mathcal{F}\subset \varphi$ the map $g_\mathcal{F}$ is
quasiminimal and contains $f$. Therefore $\mathcal{F}\in\lambda(f)$.

For the inverse implication assume that $\mathcal{F}\in\lambda(f)$. It is
easy to see that
\begin{equation}\label{tOrderConvEq1}
\underline{g}_\mathcal{F}(x)=\sup\{\underline{g}(x):g \mbox{ is a
usco map}, [g]\in\mathcal{F}\}.
\end{equation}%
For a given usco map $g$, denote by $\alpha_g$ the unique minimal
usco contained in the map $x\mapsto
[\underline{g}(x),\underline{g}^*(x)]$, where $\underline{g}^*$ is
the upper semicontinuous envelope of $\underline{g}$, i.e.
$\underline{g}^\ast$ is the pointwise infimum of all upper
semicontinuous functions greater then $\underline{g}$. Clearly
$\alpha_g$ is a lower bound of $[g]$. Hence the set $\mathcal{F}^-$
is not empty. Furthermore $f=\alpha_{g_\mathcal{F}}$. Then it
follows from (\ref{tOrderConvEq1}) that
\begin{equation}\label{tOrderConvEq2}
f=\sup\{\alpha_g:g \mbox{ is a usco map},[g]\in\mathcal{F}\}\le\sup\mathcal{F}^-.
\end{equation}
In a similar way we prove that $\mathcal{F}^+$ is not empty and
that
\begin{equation}\label{tOrderConvEq3}
f\geq \inf \mathcal{F}^+.
\end{equation}
Using that $\sup \mathcal{F}^-\leq \inf \mathcal{F}^+$ and the
inequalities (\ref{tOrderConvEq2}) and (\ref{tOrderConvEq3}) we
obtain (\ref{f=sup=inf}). Hence $\mathcal{F}$ order converges to
$f$.
\end{proof}

\begin{remark}\label{remPointInfSup}
Let us note that the infimum and the supremum in (\ref{f=sup=inf})
are not the pointwise ones. More precisely, we have that $f$ is
the unique minimal usco map contained in the quasiminimal usco map
\[
x\mapsto[\sup_{\phi\in\mathcal{F}^-}\underline{\phi}(x),\inf_{\psi\in\mathcal{F}^+}\overline{\psi}(x)]\
,\ \ x\in X.
\]
Furthermore, if $X$ is a Baire space there exists a residual
subset $D$ of $X$ such that for all $x\in D$ the value of $f$ is a
singleton and
\[
f(x)=\sup_{\phi\in\mathcal{F}^-}\underline{\phi}(x)=\inf_{\psi\in\mathcal{F}^+}\overline{\psi}(x)
\]
\end{remark}

\begin{remark}\label{remOrderConv1}
The concept of order convergence is better known in the context of
sequences, \cite{RieszI}. Let us recall that a sequence $(f_n)$ on
$\mathcal{M}(X,\mathbb{R})$ order converges to
$f\in\mathcal{M}(X,\mathbb{R})$ if there exist an increasing
sequence $(\alpha_n)$ and a decreasing sequence $(\beta_n)$ on
$\mathcal{M}(X,\mathbb{R})$ such that
\begin{eqnarray}
&&\alpha_n\leq f_n\leq \beta_n,\label{OrderConvCond1}\\
&&\sup\alpha_n=\inf\beta_n=f.\label{OrderConvCond2}
\end{eqnarray}
Using the Dedekind completeness of $\mathcal{M}(X,\mathbb{R})$ it
is easy to see that the order convergence of filters given through
(\ref{f=sup=inf}) induces the order convergence of sequences
defined above. Therefore, the class of order convergent sequences
coincides with the class of convergence sequences in $\lambda$.
\end{remark}
\begin{remark}
It was shown in \cite{AnguelovWalt} that the sequential order
convergence on $\mathbb{H}(X,\mathbb{R})$ cannot be induced by
topology. Using that the mapping $f\mapsto f_C$ is an order
isomorphism from $\mathcal{M}(X,\mathbb{R})$ to
$\mathbb{H}(X,\mathbb{R})$, this also holds true for the order
convergence on $\mathcal{M}(X,\mathbb{R})$. Since the convergence
structure $\lambda$ induces the sequential order convergence on
$\mathcal{M}(X,\mathbb{R})$, see Theorem \ref{tOrderConv} and
Remark \ref{remOrderConv1}, the convergence structure $\lambda$ on
$\mathcal{M}(X,\mathbb{R})$ is not topological.
\end{remark}

\section{Uniform convergence structure on
$\mathcal{M}(X,Y)$.}\label{sec-unif}

In this section we assume that $X$ is a Baire space and $Y$ is a
metric space with a metric $\rho$. In this case a usco mapping
$f:X\to Y$ is quasiminimal if and only if it is singlevalued at
points of a dense (equivalently residual) set
(Lemma~\ref{tQminDense} and~\ref{q-uniq}). We will need the
following lemma on product mappings.

\begin{lemma}\label{product} Let $f$ and $g$ be usco mappings from
$X$ to $Y$ and $f\times g:X\to Y\times Y$ be defined by $(f\times
g)(x)=f(x)\times g(x)$. Then the following is true.
\begin{itemize}
    \item[(i)] $f\times g$ is usco.
    \item[(ii)] If $f$ and $g$ are quasiminimal then $f\times g$ is
    quasiminimal as well.
    \item[(iii)] $f\times g$ need not be minimal even if $f$ and $g$
    are minimal.
\end{itemize}
\end{lemma}

\begin{proof} The assertion (i) is well-known and easy to see. To
show (ii) it is enough to notice that $f\times g$ is singlevalued at
points of a residual set whenever both $f$ and $g$ have that
property.

To show (iii) set $X=[0,\omega]$, $Y=[0,1]$ and define
$$f(x)=g(x)=\begin{cases} \{0\}, & x<\omega\mbox{ even},\\ \{1\}, &
x<\omega\mbox{ odd}, \\ \{0,1\}, & x=\omega. \end{cases}$$
Then $f$ and $g$ are minimal but $f\times g$ is not minimal.
\end{proof}

In particular, if $f,g\in\mathcal{M}(X,Y)$ the product mapping
$f\times g$ is quasiminimal. So we can define a mapping
$\chi:\mathcal{M}(X,Y)\times\mathcal{M}(X,Y)\to\mathcal{M}(X,Y\times
Y)$ by the formula
$$\{\chi(f,g)\}=[f\times g].$$

Further, if $\phi:X\to Y\times Y$ is usco, the composed mapping
$\rho\circ\phi:X\to\mathbb{R}$ defined by
$$(\rho\circ\phi)(x)=\rho(\phi(x))=\{\rho(y_1,y_2):(y_1,y_2)\in
\phi(x)\}$$
is usco as well. This follows from the fact that the metric
$\rho:Y\times Y\to\mathbb{R}$ is a continuous mapping.

Now we are ready to define a uniform convergence structure on
$\mathcal{M}(X,Y)$. Let us recall that such a uniform convergence
structure is a collection $\Upsilon$ of filters on
$\mathcal{M}(X,Y)\times\mathcal{M}(X,Y)$ satisfying the following
conditions (see \cite{Beattie}):
\begin{eqnarray}
\bullet&&\langle f\rangle\times\langle f\rangle\in\Upsilon\ \mbox{
for all }\
f\in\mathcal{M}(X,Y).\label{UniformCond1}\\
\bullet&&\mathcal{U}\cap\mathcal{V}\in\Upsilon\ \mbox{ whenever }\
\mathcal{U},\mathcal{V}\in\Upsilon.\label{UniformCond2}\\
\bullet&&\mbox{If }\ \mathcal{U}\in\Upsilon \ ,\  \mbox{ then }\
\mathcal{V}\in\Upsilon\ \mbox{ for each filter $\mathcal{V}$ on }\nonumber\\
&&\mathcal{M}(X,Y)\times\mathcal{M}(X,Y)\ \mbox{ such that }\
\mathcal{V}\supseteq\mathcal{U}.\label{UniformCond3}\\
\bullet&&\mbox{If }\ \mathcal{U}\in\Upsilon\ \mbox{ then }\
\mathcal{U}^{-1}\in\Upsilon.\label{UniformCond4}\\
\bullet&&\mbox{For all } \mathcal{U},\mathcal{V}\in\Upsilon\ \mbox{
one has }\ \mathcal{U}\circ\mathcal{V}\in\Upsilon\ \mbox{
whenever}\nonumber\\
&&\mbox{the composition }\ \mathcal{U}\circ\mathcal{V}\ \mbox{
exists}.\label{UniformCond5}
\end{eqnarray}

Recall that $\langle f\rangle$ denotes the filter generated by
$\{\{f\}\}$ and that if $\mathcal F_1$ and $\mathcal F_2$ are
filters on $\mathcal{M}(X,Y)$, $\mathcal F_1\times\mathcal F_2$
denotes the filter on $\mathcal{M}(X,Y)\times\mathcal{M}(X,Y)$ which
is generated by $\{F_1\times F_2: F_1\in \mathcal{F}_1, F_2\in
\mathcal{F}_2\}$.

 In
(\ref{UniformCond4}) above we use the common notation: If $U$ is a
subset of $\mathcal{M}(X,Y)\times\mathcal{M}(X,Y)$ then
\[
U^{-1}=\{(f,g):(g,f)\in U\}.
\]
For any filter $\mathcal{U}$ on
$\mathcal{M}(X,Y)\times\mathcal{M}(X,Y)$ we have
$\mathcal{U}^{-1}=\{U^{-1}:U\in\mathcal{U}\}$. The operation
composition used in (\ref{UniformCond5}) is defined as follows. For
any two subsets $U$ and $V$ of
$\mathcal{M}(X,Y)\times\mathcal{M}(X,Y)$
\begin{multline*}
U\circ V=\{(f,g)\in\mathcal{M}(X,Y)\times\mathcal{M}(X,Y): \\
\exists h\in\mathcal{M}(X,Y):(f,h)\in V,\, (h,g)\in U\}.
\end{multline*}
If $\mathcal{U}$ and $\mathcal{V}$ are filters on
$\mathcal{M}(X,Y)\times\mathcal{M}(X,Y)$ and $U\circ V\neq\emptyset$
for all $U\in\mathcal{U}$ and all $V\in\mathcal{V}$ then the filter
generated by $\{U\circ V:U\in\mathcal{U},\, V\in\mathcal{V}\}$ is
denoted by $\mathcal{U}\circ\mathcal{V}$ and called the composition
filter of $\mathcal{U}$ and $\mathcal{V}$. In this case one says
that the composition $\mathcal{U}\circ\mathcal{V}$ exists.

 Let $\Upsilon$ be the family of all filters
$\mathcal{U}$ on $\mathcal{M}(X,Y)\times\mathcal{M}(X,Y)$ such that
\begin{itemize}
    \item The filter $\chi(\mathcal{U})$ is bounded in
    $\mathcal{M}(X,Y\times Y)$.
    \item The filter generated by the family
    $$\{[\rho\circ\phi]:\phi\mbox{ is
    usco},[\phi]\in\chi(\mathcal{U})\}$$
    converges to $0$ (i.e. to the constant function equal to $0$) in
    $\mathcal{M}(X,\mathbb{R})$.
\end{itemize}

Let us remark that by  $\chi(\mathcal{U})$ we denote the filter
generated by $\{\chi(U):U\in\mathcal{U}\}$. We will show that
$\Upsilon$ is a uniform convergence structure on $\mathcal{M}(X,Y)$
which induces the convergence structure defined in
Definition~\ref{defFilterConv}. To do this we will need some
lemmata.

First we give an equivalent description of $\Upsilon$. Denote by
$\Delta$ the diagonal in $Y\times Y$ and set
$$\mathcal{D}=\{\phi\in\mathcal{M}(X,Y\times
Y):\phi(X)\subset\Delta\}.$$
Then we have the following.

\begin{lemma}\label{diagonal}
Let $\mathcal{U}$ be a filter on
$\mathcal{M}(X,Y)\times\mathcal{M}(X,Y)$ such that the filter
$\chi(\mathcal{U})$ is bounded in $\mathcal{M}(X,Y\times Y)$. Then
the following are equivalent.
\begin{enumerate}
    \item[(a)] $\mathcal{U}\in\Upsilon$.
    \item[(b)] $[g_{\chi(\mathcal{U})}]\subset\mathcal{D}$.
\end{enumerate}
\end{lemma}

\begin{proof} First let us show the equality
\begin{equation}\label{ppp}
\bigcap\{\rho\circ\phi: \phi\mbox{ is a
usco},[\phi]\in\chi(\mathcal{U})\}= \rho\circ\bigcap\{\phi:
\phi\mbox{ is a usco},[\phi]\in\chi(\mathcal{U})\}.
\end{equation}
The inclusion $\supset$ is obvious. Let us show the inverse one. Let
$x\in X$ and $$t\in\bigcap\{\rho(\phi(x)): \phi\mbox{ is a
usco},[\phi]\in\chi(\mathcal{U})\}.$$ For each $\phi$ such that
$[\phi]\in\chi(\mathcal{U})$ there are some
$(y^1_\phi,y^2_\phi)\in\phi(x)$ such that
$t=\rho(y^1_\phi,y^2_\phi)$. Consider the $\phi$'s ordered by the
inverse inclusion. Then the pairs $(y^1_\phi,y^2_\phi)$ form a net.
As $(y^1_\phi,y^2_\phi)\in\psi(x)$ for every $\phi\subset\psi$ and
$\psi(x)$ is compact, we can without loss of generality suppose that
the net converges to some $(y^1,y^2)\in Y\times Y$. Then clearly
$t=\rho(y^1,y^2)$ and, moreover $$(y^1,y^2)\in\bigcap\{\phi(x):
\phi\mbox{ is a usco},[\phi]\in\chi(\mathcal{U})\},$$ hence
$$t\in\rho\left(\bigcap\{\phi(x): \phi\mbox{ is a
usco},[\phi]\in\chi(\mathcal{U})\}\right),$$ which we wanted to
prove.

Let us proceed to the proof of the equivalence of (a) and (b). Let
(a) hold. It means that the usco on the left-hand side of
(\ref{ppp}) is quasiminimal and the unique minimal usco contained in
it is the constant zero function. By (\ref{ppp}) the same is true
for the usco on the right-hand side which is equal to $\rho\circ
g_{\chi(\mathcal{U})}$. We will show that (b) is satisfied. Suppose
that $\sigma\in[g_{\chi(\mathcal{U})}]\setminus\mathcal{D}$. It
means that there is $x_0\in X$ and distinct points $y_1,y_2\in Y$
such that $(y_1,y_2)\in\sigma(x_0)$. Choose disjoint open sets $U_1$
and $U_2$ in $Y$ such that $y_i\in U_i$ for $i=1,2$. As $\sigma$ is
minimal there is, by Lemma~\ref{minimal-prop}, a nonempty open set
$V\subset X$ such that $\sigma(V)\subset U_1\times U_2$. Hence
$0\notin\rho(\sigma(x))$ for $x\in V$. Therefore,
$[\rho\circ\sigma]$ is a nonempty subset of $[\rho\circ
g_{\chi(\mathcal{U})}]$ not containing the constant zero function.
It is a contradiction.

Conversely suppose that (b) holds. By (\ref{ppp}) it is enough to
prove that $[\rho\circ g_{\chi(\mathcal{U})}]$ contains only the
zero function. Suppose that $\alpha$ is a non-zero element of this
set. Then, by the minimality of $\alpha$ and
Lemma~\ref{minimal-prop} there is $c>0$ and nonempty open set
$V\subset X$ such that $\alpha(V)\subset (c,+\infty)$. Let
$$H=\{(y_1,y_2)\in Y\times Y:\rho(y_1,y_2)\ge c\}.$$
Then $H$ is a closed subset of $Y\times Y$ and
$g_{\chi(\mathcal{U})}(x)\cap H\ne\emptyset$ for each $x\in V$.
Therefore the mapping
$$h(x)=\begin{cases}g_{\chi(\mathcal{U})}(x)\cap H,& x\in V,\\
g_{\chi(\mathcal{U})}(x),& x\in X\setminus V,\end{cases}$$
is a usco map (Lemma~\ref{usco-modif}) contained in
$g_{\chi(\mathcal{U})}$. However, $[h]$ does not intersect
$\mathcal{D}$, a contradiction.
\end{proof}

\begin{lemma}\label{G-product}
Let $\mathcal{F}_1$ and $\mathcal{F}_2$ be filters on
$\mathcal{M}(X,Y)$.
\begin{itemize}
    \item [(i)] The filter $\chi(\mathcal{F}_1\times\mathcal{F}_2)$
    is bounded in $\mathcal{M}(X,Y\times Y)$ if and only if both $\mathcal{F}_1$ and $\mathcal{F}_2$
    are bounded in $\mathcal{M}(X,Y)$.
    \item [(ii)]
    $\mathcal{G}_{\chi(\mathcal{F}_1\times\mathcal{F}_2)} \subset\chi(\mathcal{G}_{\mathcal{F}_1}
    \times\mathcal{G}_{\mathcal{F}_2})$.
\end{itemize}
\end{lemma}

\begin{proof}
(i) Suppose that $\mathcal{F}_1$ and $\mathcal{F}_2$
    are bounded in $\mathcal{M}(X,Y)$. Then there are usco maps
    $g_1$ and $g_2$ such that $[g_1]\in\mathcal{F}_1$ and
    $[g_2]\in\mathcal{F}_2$. Then
    $$[g_1\times
    g_2]\supset\chi([g_1]\times[g_2])\in\chi(\mathcal{F}_1\times\mathcal{F}_2),$$
    hence $\chi(\mathcal{F}_1\times\mathcal{F}_2)$ is bounded, too.

    Conversely, let $\chi(\mathcal{F}_1\times\mathcal{F}_2)$ be
    bounded. Then there is a usco map $\phi$ with
    $[\phi]\in\chi(\mathcal{F}_1\times\mathcal{F}_2)$. Denote by
    $p_1$ and $p_2$ the projections of $Y\times Y$ onto the first
    and second coordinates, respectively. Then $g_j=p_j\circ\phi$ is
    a usco mapping from $X$ to $Y$ for $j=1,2$. Moreover,
    $\phi\subset g_1\times g_2$. Hence $[g_1\times
    g_2]\in\chi(\mathcal{F}_1\times\mathcal{F}_2)$. It means that
    there are $F_j\in\mathcal{F}_j$  for $j=1,2$ such that
    $\chi(F_1\times F_2)\subset [g_1\times g_2]$.
    We claim that $F_1\times F_2\subset[g_1]\times [g_2]$. Let $(h_1,h_2)\in F_1\times F_2$.
    By Lemma~\ref{q-uniq} there is a dense subset $D$ of $X$ such that for each $x\in
    D$ both
    $h_1(x)$ and $h_2(x)$ are singletons. Hence, for each $x\in D$ we have
    $$h_1(x)\times h_2(x)=\chi(h_1,h_2)(x)\subset g_1(x)\times
    g_2(x),$$
    so $h_j(x)\subset g_j(x)$
    for $j=1,2$. It follows from Lemma~\ref{tfgequal} that $f_j\subset g_j$ for $j=1,2$.

    Hence, for $i=1,2$ we have $F_j\subset[g_j]$. Denote by $f_j$ the set-valued mapping obtained
    as the closure of $\bigcup F_j$ in $X\times Y$. By \cite[Lemma 3.1.1]{fabian} it is usco.
     Now, clearly  $F_j\subset[f_j]$ and so $\mathcal{F}_j$ is bounded.

 (ii) Let $U\in\mathcal{G}_{\chi(\mathcal{F}_1\times\mathcal{F}_2)}$. Then there is
     a usco mapping $\phi$ with
     $[\phi]\in\chi(\mathcal{F}_1\times\mathcal{F}_2)$ such that
     $[\phi]\subset U$. Further, there are $F_j\in\mathcal{F}_j$ for
     $j=1,2$, such that $\chi(F_1\times F_2)\subset[\phi]$.

     Denote by $f_j$ the set-valued mapping obtained as the closure
     of $\bigcup F_j$ in $X\times Y$. In the same way as in the
     proof of (i) we can show that $f_j$ is a usco map. As
     $F_j\subset[f_j]$, we have $[f_j]\in\mathcal{F}_j$. Therefore we will be done if we show that
     $$\chi([f_1]\times[f_2])\subset[\phi].$$
     Let $h\in\chi([f_1]\times[f_2])$. Then $h$ is a minimal usco
     and $h\subset f_1\times f_2$. Suppose that $h\not\subset\phi$.
     Choose $x_0\in X$ and $(y_0,z_0)\in h(x_0)\setminus \phi(x_0)$. Let
     $V_1$ and $V_2$ be disjoint open subset of $Y\times Y$ with
     $\phi(x_0)\subset V_1$ and $(y_0,z_0)\in V_2$. Choose $W_1$,
     $W_2$ open subsets in $Y$ such that $(y_0,z_0)\in W_1\times
     W_2\subset V_2$.

     As $\phi$ is usco, there is $U_0$, a neighborhood of $x_0$ such
     that $\phi(U_0)\subset V_1$. As $h$ is minimal, there is (by
     Lemma~\ref{minimal-prop}) a nonempty open set $U_1\subset U_0$
     with $h(U_1)\subset W_1\times W_2$. Choose some
     $x_1\in U_1$ and $(y_1,z_1)\in h(x_1)$. Then $y_1\in f_1(x_1)$.
     By the definition of $f_1$ there is some $g_1\in F_1$, $x_2\in
     U_1$ and $y_2\in g_1(x_2)\cap W_1$. As $g_1$ is minimal, there
     is (again by Lemma~\ref{minimal-prop}) a nonempty open set
     $U_2\subset U_1$ with $g_1(U_2)\subset W_1$.
     Similarly there is some $g_2\in F_2$ and a nonempty open set
     $U_3\subset U_2$ such that $g_2(U_3)\subset W_2$.
     Thus $(g_1\times g_2)(U_3)\subset W_1\times W_2$, so
     $(g_1\times g_2)(U_3)\cap \phi(U_3)=\emptyset$. Therefore
     $\chi(g_1,g_2)\notin[\phi]$, a contradiction.
\end{proof}

\begin{lemma}\label{induced} Let $f\in\mathcal{M}(X,Y)$ and
$\mathcal{F}$ be a filter on $\mathcal{M}(X,Y)$. Then $\langle
f\rangle\times\mathcal{F}\in\Upsilon$ if and only if $\mathcal{F}\to
f$ in $\mathcal{M}(X,Y)$.
\end{lemma}

\begin{proof}
Suppose that $\mathcal{F}\to f$. Then $\mathcal{F}$ is bounded.
Moreover, $\langle f\rangle$ is also bounded, hence $\chi(\langle
f\rangle\times\mathcal{F})$ is bounded as well by
Lemma~\ref{G-product}.

Moreover, if $g$ is a usco map such that $[g]\in\mathcal{F}$, then
$[f\times g]\in\chi(\langle f\rangle\times\mathcal{F})$. Thus
$$g_{\chi(\langle
f\rangle\times\mathcal{F})} \subset \bigcap\{f\times g: g\mbox{ is
usco}, [g]\in\mathcal{F}\} = f\times g_{\mathcal{F}}.$$
As $[g_{\mathcal{F}}]=\{f\}$, we have
$$[g_{\chi(\langle
f\rangle\times\mathcal{F})}]=[f\times f].$$
As the diagonal $\Delta$ is closed in $Y\times Y$ and $(f\times
f)(x)\cap\Delta\ne\emptyset$ for each $x\in D$, the mapping
$x\mapsto (f\times f)(x)\cap\Delta$ is usco (by
Lemma~\ref{usco-modif}). Hence $\chi(f,f)(x)\subset \Delta$ for each
$x\in D$, i.e. $\chi(f,f)\in\mathcal{D}$. By Lemma~\ref{diagonal}
this completes the proof that $\chi(\langle
f\rangle\times\mathcal{F})$ belongs to $\Upsilon$.

Conversely, suppose that $\langle f\rangle\times\mathcal{F}$ belongs
to $\Upsilon$. Then the filter $\chi(\langle
f\rangle\times\mathcal{F})$ is bounded, and hence $\mathcal{F}$ is
bounded as well by Lemma~\ref{G-product}. By  Lemma~\ref{diagonal}
we have
$$[g_{\chi(\langle f\rangle\times\mathcal{F})}]\subset\mathcal{D}.$$
Moreover, by Lemma~\ref{G-product}(ii) we get
\begin{align*}
[g_{\chi(\langle f\rangle\times\mathcal{F})}]&\supset
\bigcap\{\chi(\{f\}\times[g]): g\mbox{ is usco},[g]\in\mathcal{F}\}
\\& \supset \chi\left(\bigcap\{\{f\}\times[g]: g\mbox{ is
usco},[g]\in\mathcal{F}\}\right)
=\chi(\{f\}\times[g_{\mathcal{F}}]),
\end{align*}
hence
$$\chi(\{f\}\times[g_{\mathcal{F}}])\subset\mathcal{D}.$$
If $h\in[g_{\mathcal{F}}]$ is different from $f$, then $f(x)\cap
h(x)=\emptyset$ for some $x\in X$ (by Lemma~\ref{tfgequal}). But
this implies that $\chi(f,h)\notin \mathcal{D}$, a contradiction.
Hence $[g_{\mathcal{F}}]=\{f\}$, i.e. $\mathcal F\to f$.
\end{proof}

\begin{theorem}\label{tUniformConv}
The collection of filters $\Upsilon$ is a uniform convergence
structure inducing the convergence structure on $\mathcal{M}(X,Y)$.
\end{theorem}

\begin{proof}
To prove that $\Upsilon$ is a uniform convergence structure we need
to show that $\Upsilon$ satisfies the properties
(\ref{UniformCond1})--(\ref{UniformCond5}). The fact that $\Upsilon$
generates the convergence structure on $\mathcal{M}(X,Y)$ then
follows immediately from Lemma~\ref{induced}.

The property (\ref{UniformCond3}) is obvious. To show the property
(\ref{UniformCond4}) it is enough to use the symmetry of the metric
$\rho$. The property (\ref{UniformCond1}) follows immediately from
Lemma~\ref{induced} as $\langle f\rangle\to f$.

Let us show the property (\ref{UniformCond2}). Let $\mathcal{U}$ and
$\mathcal{V}$ belong to $\Upsilon$. First we show that
$$\chi(\mathcal{U}\cap\mathcal{V})=\chi(\mathcal{U})\cap\chi(\mathcal{V}).$$
Indeed, the inclusion $\subset$ is obvious. To prove the inverse one
choose an element $S$ in the set on the right-hand side. Then there
are $U\in \mathcal{U}$ and $V\in\mathcal{V}$ such that
$\chi(U)\subset S$ and $\chi(V)\subset S$. Then $U\cup
V\in\mathcal{U}\cap\mathcal{V}$ and
$$\chi(U\cup V)=\chi(U)\cup\chi(V)\subset S,$$
hence $S\in\chi(\mathcal{U}\cap\mathcal{V})$. Now, in the same way
as in the proof of Theorem~\ref{tConvSpace} one can easily show that
$\chi(\mathcal{U}\cap\mathcal{V})$ is bounded and, moreover,
$$g_{\chi(\mathcal{U}\cap\mathcal{V})}\subset g_{\chi(\mathcal{U})}\cup
g_{\chi(\mathcal{V})}.$$
Therefore, by Lemma~\ref{diagonal} it is enough to prove the
following claim:
\begin{equation}\label{d-union}
\phi,\psi:X\to Y\times Y \mbox{ usco maps},
[\phi]\cup[\psi]\subset\mathcal{D}\Rightarrow[\phi\cup\psi]\subset\mathcal{D}.
\end{equation}
Suppose that $h\in[\phi\cup\psi]\setminus \mathcal{D}$. Choose
$x_0\in X$ and $(y_1,y_2)\in h(x_0)$ such that $y_1\ne y_2$. Find
disjoint open sets $V_1,V_2\subset Y$ such that $y_1\in V_1$ and
$y_2\in V_2$. By Lemma~\ref{minimal-prop} there is a nonempty open
set $U_0\subset X$ such that $h(U_0)\subset V_1\times V_2$.

We claim that there is a nonempty open set $U_1\subset U_0$ such
that either $h|_{U_1}\subset\phi|_{U_1}$ or
$h|_{U_1}\subset\psi|_{U_1}$. Indeed, suppose that
$h|_{U_0}\not\subset\phi|_{U_0}$. As $h|_{U_0}$ is minimal
(\cite[Lemma 2]{stenfra}) by Lemma~\ref{tfgequal} we get a nonempty
open set $U_1\subset U_0$ with $\phi(U_1)\cap h(U_1)=\emptyset$. As
$h\subset\phi\cup\psi$, it follows $h|_{U_1}\subset\psi|_{U_1}$
which proves our claim.

So suppose, say, that $h|_{U_1}\subset\phi|_{U_1}$. Define a mapping
$\widetilde\phi$ by
$$\widetilde\phi(x)=\begin{cases}\phi(x),& x\in X\setminus U_1,\\
h(x),& x\in U_1.\end{cases}$$
By Lemma~\ref{usco-modif} it is a usco map. Further,
$[\widetilde\phi]\subset[\phi]$ and
$[\widetilde\phi]\cap\mathcal{D}=\emptyset$ (as
$\widetilde\phi(U_1)\subset V_1\times V_2$), a contradiction
completing the proof of (\ref{d-union}).

It remains to prove the condition (\ref{UniformCond5}). Let
$\mathcal{U}$ and $\mathcal{V}$ be elements of $\Upsilon$ such that
$\mathcal{U}\circ\mathcal{V}$ exists.

First let us show that $\chi(\mathcal{U}\circ\mathcal{V})$ is
bounded. We know that both $\chi(\mathcal{U})$ and
$\chi(\mathcal{V})$ are bounded, and hence
$\chi(\mathcal{U}\cap\mathcal{V})$ is bounded as well (by the
already proved condition (\ref{UniformCond2})). Hence there is a
usco map $\phi$ such that
$[\phi]\in\chi(\mathcal{U}\cap\mathcal{V})$. Let
$\alpha=p_1\circ\phi$ and $\alpha=p_2\circ\phi$ (where $p_1$ and
$p_2$ are projections of $Y\times Y$, see the proof of
Lemma~\ref{G-product}(i)). Then $\alpha$ and $\beta$ are usco maps
and $[\alpha\times\beta]\in\chi(\mathcal{U}\cap\mathcal{V})$. Hence
 there is some $U\in\mathcal{U}\cap\mathcal{V}$ such that
$\chi(U)\subset[\alpha\times\beta]$. We will show that $\chi(U\circ
U)\subset[\alpha\times\beta]$ as well.

Let $(f,g)\in U\circ U$. Then there is $h\in\mathcal{M}(X,Y)$ such
that $(f,h)\in U$ and $(h,g)\in U$. Thus both $\chi(f,h)$ and
$\chi(h,g)$ are contained in $\alpha\times \beta$. By
Lemma~\ref{q-uniq} there is a dense set $D\subset X$ such that for
each $x\in D$ all the values $f(x)$, $h(x)$ and $g(x)$ are
singletons. Hence for $x\in D$ we have $f(x)\subset\alpha(x)$ and
$g(x)\subset\beta(x)$. In particular, $\chi(f,g)(x)=f(x)\times
g(x)\subset(\alpha\times\beta)(x)$ for $x\in D$. Therefore
$\chi(f,g)\subset \alpha\times\beta$ by Lemma~\ref{tfgequal}.

This completes the proof that $\chi(U\circ
U)\subset[\alpha\times\beta]$ and hence
$\chi(\mathcal{U}\circ\mathcal{V})$ is bounded.

To finish the proof that $\mathcal{U}\circ\mathcal{V}$ belongs to
$\Upsilon$ we will use Lemma~\ref{diagonal}. Suppose that
$\alpha\in\mathcal{M}(X,Y\times Y)\setminus\mathcal{D}$. Choose
$x_0\in X$ and distinct points $y_0,z_0\in Y$  such that
$(y_0,z_0)\in \alpha(x_0)$. Let $c>0$ be such that
$c<\rho(y_0,z_0)$. By the minimality of $\alpha$ and
Lemma~\ref{minimal-prop} there is a nonempty open set $U_0\subset X$
such that
$$\alpha(U_0)\subset \{(y,z)\in Y\times Y:\rho(y,z)>c\}.$$
By the already proved condition (\ref{UniformCond2}) we know that
$[g_{\chi(\mathcal{U}\cap\mathcal{V})}]\subset\mathcal{D}$. Thus
there is some $x_1\in U_0$ such that
$$g_{\chi(\mathcal{U}\cap\mathcal{V})}(x_1)\cap
\{(y,z)\in Y\times Y:\rho(y,z)\ge\tfrac c2\}=\emptyset.$$
Indeed, otherwise
$$h(x)=\begin{cases}g_{\chi(\mathcal{U}\cap\mathcal{V})}(x)\cap \{(y,z)\in Y\times Y:\rho(y,z)\ge\frac c2\}, &
x\in U_0, \\ g_{\chi(\mathcal{U}\cap\mathcal{V})}(x), & x\in
X\setminus U_0,\end{cases}$$
would be a usco mapping (by Lemma~\ref{usco-modif}) contained in
$g_{\chi(\mathcal{U}\cap\mathcal{V})}$ but not containing any
element of $\mathcal{D}$, a contradiction. Now, by the definition of
$g_{\chi(\mathcal{U}\cap\mathcal{V})}$ there is some usco map $\phi$
with $[\phi]\in\chi(\mathcal{U}\cap\mathcal{V})$ such that
$$\phi(x_1)\cap\{(y,z)\in Y\times Y:\rho(y,z)\ge\tfrac
c2\}=\emptyset.$$ As $\phi$ is usco, there is an open set $U_1$ with
$x_1\in U_1\subset U_0$ such that $$\phi(U_1)\cap\{(y,z)\in Y\times
Y:\rho(y,z)\ge\tfrac c2\}=\emptyset.$$

There is some $M\in\mathcal{U}\cap\mathcal{V}$ such that
$\chi(M)\subset [\phi]$. Then we have
\begin{equation}\label{star} \chi(M\circ
M)\subset[\phi\star\phi],\end{equation}
where $\phi\star\phi$ is the usco mapping defined by
$$(\phi\star\phi)(x)=\phi(x)\circ\phi(x).$$
Let us show first that $\phi\star\phi$ is a usco mapping. We will
use \cite[Lemma 3.1.1]{fabian}. Let $x_\tau$ be a net in $X$
converging to some $x\in X$ and let
$(y_\tau,z_\tau)\in(\phi\star\phi)(x_\tau)$. For each $\tau$ there
is some $u_\tau\in Y$ such that $(y_\tau,u_\tau)\in\phi(x_\tau)$ and
$(u_\tau,z_\tau)\in\phi(x_\tau)$. As $\phi$ is usco, there is a
subnet $(y_\nu,u_\nu)$ of $(y_\tau,u_\tau)$ converging to some
$(y,u)\in\phi(x)$. Using once more that $\phi$ is usco, we obtain a
subnet $(u_\mu,z_\mu)$ of $(u_\nu,z_\nu)$ converging to some
$(u,z)\in\phi(x)$ (note that $u_\nu$ converges to $u$, and hence the
limit of $u_\mu$ is also $u$). Then $(x_\mu,z_\mu)$ converges to
$(x,z)$ and $(x,z)\in\phi(x)\circ\phi(x)$.

Let us proceed to the proof of (\ref{star}). Pick $(f,g)\in M\circ
M$. Then there is $h\in\mathcal{M}(X,Y)$ such that both $(f,h)$ and
$(h,g)$ belong to $M$. Hence both $\chi(f,h)$ and $\chi(h,g)$ are
contained in $\phi$. By Lemma~\ref{q-uniq} there is a dense set
$D\subset X$ such that all the mappings $f,g,h$ are singlevalued on
$D$. Let $x\in D$. Then $f(x)\times h(x)\subset \phi(x)$ and
$h(x)\times g(x)\subset \phi(x)$, hence $$f(x)\times
g(x)\subset\phi(x)\circ\phi(x)=(\phi\star\phi)(x).$$ Therefore
$\chi(f,g)(x)\subset(\phi\star\phi)(x)$ for each $x\in D$. It
follows from Lemma~\ref{tfgequal} that
$\chi(f,g)\subset\phi\star\phi$ which completes the proof of
(\ref{star}).

Hence $\phi\star\phi$ is a usco map and
$[\phi\star\phi]\in\chi(\mathcal{U}\circ\mathcal{V})$. Let $x\in
U_1$ and $(y,z)\in(\phi\star\phi)(x)$. Then there is $u\in Y$ such
that both $(y,u)$ and $(u,z)$ belong to $\phi(x)$. Then
$$\rho(y,z)\le\rho(y,u)+\rho(u,z)<c.$$
Thus $(\phi\star\phi)(U_1)\cap\alpha(U_1)=\emptyset$. Therefore
$\alpha\notin[\phi\star\phi]$ which completes the proof.
\end{proof}

An important question associated with uniform convergence spaces
is their completeness, that is, the convergence of Cauchy filters.
Let us recall that a filter $\mathcal{F}$ on $\mathcal{M}(X,Y)$ is
called Cauchy if $\mathcal{F}\times\mathcal{F}\in\Upsilon$.

\begin{theorem}\label{tcomplete}
The uniform convergence space $(\mathcal{M}(X,Y),\Upsilon)$ is
complete.
\end{theorem}
\begin{proof}
Assume that the filter $\mathcal{F}$ on $\mathcal{M}(X,Y)$  is
Cauchy, that is, $\mathcal{F}\times\mathcal{F}\in\Upsilon$. Then
$\mathcal{F}$ is bounded by Lemma~\ref{G-product}(i). Hence the usco
map $g_\mathcal{F}$ is well defined and nonempty valued on $X$.
Moreover, by Lemma~\ref{diagonal} and Lemma~\ref{G-product}(ii) we
have
\begin{align*}
\mathcal{D}&\supset[g_{\chi(\mathcal{F}\times\mathcal{F})}]\supset
\bigcap\{\chi([g]\times[g]): g\mbox{ is usco},[g]\in\mathcal{F}\}
\\& \supset \chi\left(\bigcap\{[g]\times[g]: g\mbox{ is
usco},[g]\in\mathcal{F}\}\right)
=\chi([g_{\mathcal{F}}]\times[g_{\mathcal{F}}]).
\end{align*}
If $f_1,f_2$ are two different elements of $[g_{\mathcal{F}}]$, by
Lemma~\ref{tfgequal} there is some $x\in X$ with $f_1(x)\cap
f_2(x)=\emptyset$, hence $\chi(f_1,f_2)\notin\mathcal{D}$. So
$g_\mathcal{F}$ is quasiminimal and hence $\mathcal{F}$ converges to
the unique element of $[g_{\mathcal{F}}]$.
\end{proof}

\begin{remarks}\quad
\begin{itemize}
    \item[1.] Although the definition of the uniform convergence structure $\Upsilon$ includes the
    metric $\rho$, it depends only on the topology of $Y$, i.e. is
    the same for all equivalent metrics on Y. This follows from
    Lemma~\ref{diagonal}.
    \item[2.] If $X$ is a singleton, then both
    $\mathcal{M}(X,Y)\times\mathcal{M}(X,Y)$ and $\mathcal{M}(X,Y\times
    Y)$ can be canonically identified with $Y\times Y$. In this case
    $\Upsilon$ consist of those filters $\mathcal{U}$ on $Y\times Y$
    such that there is a compact set $K\subset Y$ such that the
    filter generated by the neighborhoods of the diagonal in
    $K\times K$ is contained in $\mathcal{U}$. In particular, if $Y$
    is compact, $\Upsilon$ coincide with the (unique) uniformity on
    $Y$.
    \item[3.] We supposed that $X$ is a Baire space and $Y$ a metric
    space. In fact, the definition of the uniform structure can be
    done whenever $X$ is a Baire space, $Y$ a completely regular space and every minimal usco from $X$ to $Y$ is
    singlevalued at points of a residual subset of $X$. Let us outline the necessary differences.

    The collection $\Upsilon$ should be defined using the equivalent
    condition from Lemma~\ref{diagonal} which does not use the
    metric $\rho$. The only place after this lemma we have used the
    metric $\rho$ is the proof of the conditions (\ref{UniformCond4}) and (\ref{UniformCond5})
    in Theorem~\ref{tUniformConv}. However, the condition (\ref{UniformCond4}) is trivial also using the alternative definition and to prove the
    condition (\ref{UniformCond5}) we could use,
    instead of the metric $\rho$, a continuous pseudometric. Indeed,
    if $Y$ is completely regular (and Hausdorff -- which we
    automatically assume), then for any two distinct points $y_0,z_0\in
    Y$ there is a continuous function $f:Y\to\mathbb{R}$ with
    $f(y_0)\ne f(z_0)$. Then $d(y,z)=|f(y)-f(z)|$ defines a
    continuous pseudometric on $Y$ such that $d(y,z)>0$. Hence the
    proof can be carried on.

    In particular, Theorems~\ref{tUniformConv} and~\ref{tcomplete}
    remain true if $X$ is a Baire space and $Y$ a
    $\sigma$-fragmented Banach space equipped with the weak topology
    (see e.g. \cite{kendmoors})
    or a dual Banach space with equipped with the weak* topology
    which belongs to the Stegall's class (see \cite[Chapter
    3]{fabian}). These classes of Banach spaces are quite large and
    include, for example, all separable spaces and all reflexive
    spaces.
    \end{itemize}
\end{remarks}

\section{The subspace of continuous functions}\label{sec-dens}

The space of minimimal usco maps $\mathcal{M}(X,Y)$ contains a
natural subspace $\mathcal{C}(X,Y)$ consisting of continuous
functions from $X$ to $Y$. In the previous section we have shown
that the convergence space $\mathcal{M}(X,Y)$ is complete for the
natural uniform convergence structure whenever $X$ is a Baire space
and $Y$ is a metric space. Therefore the closure of
$\mathcal{C}(X,Y)$ in $\mathcal{M}(X,Y)$ could be viewed as a
completion of $\mathcal{C}(X,Y)$. In this section we study the
question when $\mathcal{C}(X,Y)$ is dense in $\mathcal{M}(X,Y)$.

Let us recall the definition of a closed subset of a convergence
space and related notions. A subset $A$ of a convergence space is
\textit{closed} if $f\in A$ whenever $\mathcal{F}$ is a filter
converging to $f$ and satisfying $A\in\mathcal{F}$. The
\textit{closure} of a set $A$ is the smallest closed set containing
$A$. And a set is \textit{dense} if its closure is the whole space.

First we note that $\mathcal{C}(X,Y)$ is not always dense in
$\mathcal{M}(X,Y)$.

\begin{example} $\mathcal{C}(\mathbb{R},\{0,1\})$ is a proper closed
subset of $\mathcal{M}(\mathbb{R},\{0,1\})$.
\end{example}

\begin{proof} The mapping $g$ defined by
$$g(x)=\begin{cases} \{0\}, & x<0, \\ \{0,1\}, & x=0, \\ \{1\}, &
x>0,\end{cases}$$
belongs to
$\mathcal{M}(\mathbb{R},\{0,1\})\setminus\mathcal{C}(\mathbb{R},\{0,1\})$.

Further, let us show that $\mathcal{C}(\mathbb{R},\{0,1\})$ is
closed. Let $\mathcal{F}$ be a filter on
$\mathcal{M}(\mathbb{R},\{0,1\})$ converging to some
$f\in\mathcal{M}(\mathbb{R},\{0,1\})$ satisfying
$\mathcal{C}(\mathbb{R},\{0,1\})\in\mathcal{F}$.

Let $g$ be a usco map such that $[g]\in\mathcal{F}$. Then $[g]$
contains an element of $\mathcal{C}(\mathbb{R},\{0,1\})$. As
$\mathcal{C}(\mathbb{R},\{0,1\})$ has only two elements (constant
function $0$ and constant function $1$), there is one of them
contained in $[g]$ for every $g$ satisfying $[g]\in\mathcal{F}$.
Therefore $[g_{\mathcal F}]$ contains an element of
$\mathcal{C}(\mathbb{R},\{0,1\})$, which implies
$f\in\mathcal{C}(\mathbb{R},\{0,1\})$.
\end{proof}

This example shows that in order to have $\mathcal{C}(X,Y)$ dense in
$\mathcal{M}(X,Y)$, some assumptions on $Y$ are needed. A natural
assumption of this kind is that $Y$ is a convex subset of a normed
linear space. A partial positive result is the following one.

\begin{theorem}\label{density} Let $X$ be a Baire metric space and $Y$ be a closed
convex subset of $\mathbb{R}^d$. Then $\mathcal{C}(X,Y)$ is dense in
$\mathcal{M}(X,Y)$.
\end{theorem}

\begin{proof} Let $g\in\mathcal{M}(X,Y)$. It follows for example
from \cite{HJT} that $g$ has a selection of the first Baire class,
i.e., there is a (single-valued) function $f:X\to Y$ which is of the
first Baire class (i.e., the pointwise limit of a sequence of
continuous functions) such that $f(x)\in g(x)$ for each $x\in X$.

By \cite[Theorem 3.3]{ok-b1} there is a usco map $h:X\to Y$ and a
sequence of continuous functions $f_n:X\to Y$ which pointwise
converges to $f$ and $f_n\subset h$ for each $n$.

Note that the sequence $f_n$ is bounded in $\mathcal{M}(X,Y)$. It
follows from Corollary~\ref{corseq} that the sequence $f_n$
converges to $g$ in $\mathcal{M}(X,Y)$. This completes the proof.
\end{proof}

We do not know whether the result on density is valid in more
general situations. Let us formulate some of these problems.

\begin{question} Let $X$ be a Baire metric space and $Y$ a
convex subset of a normed linear space. Is $\mathcal{C}(X,Y)$ dense
in $\mathcal{M}(X,Y)$?
\end{question}

Note that in this situation every $g\in\mathcal{M}(X,Y)$ has a
selection of the first Baire class (this follows for example from
\cite[Theorem 2.2]{srivatsa}) and Corollary~\ref{corseq} could be
applied as well. The missing ingredient is the analogue of
\cite[Theorem 3.3]{ok-b1}. It seems to be unknown whether such an
analogue holds.

Another problem is whether we can drop the assumption of
metrizability of $X$.

\begin{question} Let $X$ be a Baire topological space and $Y$ a
convex subset of a normed linear space. Is $\mathcal{C}(X,Y)$ dense
in $\mathcal{M}(X,Y)$?
\end{question}

In this case sequences are not enough as we can see from the
following example. However, to prove the density of
$\mathcal{C}(X,Y)$ we are not obliged to use sequences. Nets or
filters are allowed as well. But then some other technics should be
used, as Theorem~\ref{pointwise} (and Corollary~\ref{corseq}) is
true only for sequences (due to Example~\ref{pw-exa}).

\begin{example} There is a compact Hausdorff space $X$ and a proper
subset $A\subset\mathcal{M}(X,[0,1])$ which contains
$\mathcal{C}(X,[0,1])$ and is closed to taking limits of sequences.
Moreover, in this case $\mathcal{C}(X,[0,1])$ is dense in
$\mathcal{M}(X,[0,1])$.
\end{example}

\begin{proof} Let $X$ be the ordinal interval $[0,\omega_1]$ and
$$A=\{ g\in\mathcal{M}(X,[0,1]) : g(\omega_1) \mbox{ is a
singleton}\}.$$
Then clearly $A\supset\mathcal{C}(X,[0,1])$. Further, $A$ is a
proper subset of $\mathcal{M}(X,[0,1])$ as the mapping
$g:[0,\omega_1]\to[0,1]$ defined by
$$g(\alpha)=\begin{cases} \{0\}, & \alpha\mbox{ odd non-limit ordinal},
\\ \{1\}, & \alpha\mbox{ even non-limit ordinal},
\\ \{0,1\}, & \alpha\mbox{ limit ordinal},
\end{cases}
$$
is minimal usco and does not belong to $A$.

Next we shall show that $A$ is closed to limits of sequences. Let
$f_n$ be a sequence from $A$ converging to some
$f\in\mathcal{M}(X,[0,1])$. It follows from Lemma~\ref{q-uniq} and
Theorem~\ref{tPointConv} that there is a residual subset $D$ of $X$
such that for all $x\in D$ the values  $f(x)$ and all $f_n(x)$'s are
singletons and, moreover, $f_n(x)\to f(x)$ in the topology of
$[0,1]$. Note that the set $D$ must contain all the isolated
ordinals.

Further note, that for any $h\in A$ there is $\alpha<\omega_1$ such
that $h(x)=h(\omega_1)$ for all $x\in[\alpha,\omega_1]$. Hence, to
each $f_n$ we can associate such an $\alpha_n$. Let $\alpha$ be the
supremum of $\alpha_n$'s. Then for each isolated ordinal
$x\in[\alpha,\omega_1]$ we have
$$f(x)=\lim_{n\to\infty} f_n(x)=\lim_{n\to\infty} f_n(\omega_1).$$
Therefore $f$ assumes for each isolated $x\in[\alpha,\omega_1]$ the
same singleton value, and hence $f(\omega_1)$ is a singleton as
well. This completes the proof that $A$ is closed to taking limits
of sequences.

Now we will prove that $\mathcal{C}(X,[0,1])$ is dense in
$\mathcal{M}(X,[0,1])$.

Take any $g\in\mathcal{M}(X,[0,1])$. For each $\alpha<\omega_1$
define the mapping
$$g_\alpha(x)=\begin{cases} g(x), & x\in[0,\alpha], \\ \{0\}, & x\in
(\alpha,\omega_1].\end{cases}$$
Then each $g_\alpha$ is a minimal usco belonging to $A$. Moreover,
the net $g_\alpha$ converges to $g$. To see this we use
Theorem~\ref{tConvSeq}. We define usco maps $h_\alpha$ by the
formula
$$h_\alpha(x)=\begin{cases} g(x), & x\in[0,\alpha], \\ g(x)\cup \{0\}, & x\in
(\alpha,\omega_1].\end{cases}$$
Then $g_\alpha\subset h_\alpha$ and $h_\beta\subset h_\alpha$ for
each $\alpha\le\beta<\omega_1$. Further, the intersection of all
$h_\alpha$'s is the usco map
$$h(x)=\begin{cases} g(x), & x\in[0,\omega_1), \\ g(x)\cup \{0\}, &
x=\omega_1.\end{cases}$$
It is clear that $h$ is quasiminimal and $[h]=\{g\}$.

This shows that $A$ is dense in $\mathcal{M}([0,\alpha],[0,1])$. We
conclude by showing that $\mathcal{C}(X,[0,1])$ is dense in $A$. Let
$g\in A$. Then there is $\alpha<\omega_1$ such that
$g(x)=g(\omega_1)$ for each $x\in(\alpha,\omega_1]$. Then
$g|_{[0,\alpha]}$ belongs to $\mathcal{M}([0,\alpha],[0,1])$. As
$[0,\alpha]$ is a metrizable compact space, there is (by the proof
of Theorem~\ref{density}) a sequence of continuous functions
$f_n:[0,\alpha]\to[0,1]$ converging to $g|_{[0,\alpha]}$ in
$\mathcal{M}([0,\alpha],[0,1])$. Extend the functions $f_n$ to
functions $h_n:X\to[0,1]$ by defining $h(x)=g(\omega_1)$ for
$x>\alpha$. Then $h_n$ are continuous and clearly converge to $g$ in
$\mathcal{M}(X,[0,1])$. This completes the proof.
\end{proof}

In fact, although the assumption that the domain space $X$ is Baire
is quite natural, we do not know the answer to the following
question.

\begin{question} Is there a topological space $X$ and a convex
subset $Y$ of a normed linear space such that $\mathcal{C}(X,Y)$ is
not dense in $\mathcal{M}(X,Y)$?
\end{question}

 %The minimal usco
%generalize the continuous functions preserving some of their
%essential properties. Hence our special interest in the set $C(X,Y)$
%of continuous functions.

\end{document}